\documentclass[12pt,a4paper]{article}
\usepackage{amsmath,amsxtra,amsfonts,amssymb,latexsym,amsthm,amscd,verbatim,amstext}
\usepackage[mathscr]{eucal}

\newcommand{\dif}{\mathrm{d}}

\input cyracc.def

\newfont{\tricyr}{wncyr10 at 12pt}
\newfont{\tricyi}{wncyi10 at 12pt}
\newfont{\tricyb}{wncyb10 at 12pt}
\newfont{\Tricyr}{wncyr10 at 13.6pt}
\newfont{\Tricyi}{wncyi10 at 13.6pt}
\newfont{\Tricyb}{wncyb10 at 13.6pt}
\newfont{\tricmr}{cmr10 at 13.6pt}
\newfont{\tricmi}{cmti10 at 13.6pt}
\newfont{\tricmb}{cmb10 at 13.6pt}

\theoremstyle{plain}
\newtheorem{Th}{Theorem}

\newtheorem{Lem}{Lemma}

\theoremstyle{definition}
\newtheorem{Def}{Definition}
\newtheorem{Not}{Remark}

\begin{document}
\centerline {{\bf Well-posedness for the Navier-Stokes equations}} 
\centerline {{\bf  with datum in Sobolev-Fourier-Lorentz spaces}} 

\vskip 0.7cm
\begin{center}
D. Q. Khai, N. M. Tri
\end{center}
\begin{center}
Institute of Mathematics, Vietnam Academy of Science and Technology\\
18 Hoang Quoc Viet, 10307  Cau Giay, Hanoi, Vietnam
\end{center}
\vskip 0.5cm 
{\bf Abstract}: In this note, for $s \in \mathbb R$ and $1 \leq p, r \leq \infty$, 
we introduce and study Sobolev-Fourier-Lorentz spaces 
$\dot{H}^s_{\mathcal{L}^{p, r}}(\mathbb{R}^d)$. In the family 
spaces $\dot{H}^s_{\mathcal{L}^{p, r}}(\mathbb{R}^d)$, the critical invariant 
spaces for the Navier-Stokes equations correspond to the value $s = \frac{d}{p} - 1$. 
When the initial datum belongs to the critical spaces  $\dot{H}^{\frac{d}{p} - 1}_{\mathcal{L}^{p,r}}(\mathbb{R}^d)$ with 
$d \geq 2, 1 \leq p <\infty$, and $1 \leq r < \infty$,
we establish the existence of local mild solutions to the Cauchy problem 
for the Navier-Stokes equations in spaces  
$L^\infty([0, T]; \dot{H}^{\frac{d}{p} - 1}_{\mathcal{L}^{p, r}}(\mathbb{R}^d))$  
with arbitrary initial value, and existence of global mild solutions in spaces  
$L^\infty([0, \infty); \dot{H}^{\frac{d}{p} - 1}_{\mathcal{L}^{p, r}}(\mathbb{R}^d))$  
when the norm of the initial value in the Besov spaces 
$\dot{B}^{\frac{d}{\tilde p} - 1, \infty}_{\mathcal{L} ^{\tilde p,\infty}}(\mathbb{R}^d)$ 
is small enough, where $\tilde p$ may take some suitable values.


\footnotetext[3]{{\it e-mail address}: triminh@math.ac.vn}
\vskip 1.0cm 
\centerline{\bf{\S1. INTRODUCTION}} 
\vskip 0,5cm 
We consider the Navier-Stokes equations (NSE) in $d$ dimensions in special setting of a viscous, homogeneous, incompressible fluid which fills the entire space and is not submitted to external forces. Thus, the equations we consider are the system: 
\begin{align} 
\left\{\begin{array}{ll} \partial _tu  = \Delta u - \nabla .(u \otimes u) - \nabla p , & \\ 
{\rm div}(u) = 0, & \\
u(0, x) = u_0, 
\end{array}\right . \notag
\end{align}
which is a condensed writing for
\begin{align} 
\left\{\begin{array}{ll} 1 \leq k \leq d, \ \  \partial _tu_k  
= \Delta u_k - \sum_{l =1}^{d}\partial_l(u_lu_k) - \partial_kp , & \\ 
\sum_{l =1}^{d}\partial_lu_l = 0, & \\
1 \leq k \leq d, \ \ u_k(0, x) = u_{0k} .
\end{array}\right . \notag
\end{align}
The unknown quantities are the velocity $u(t, x)=(u_1(t, x),\dots,u_d(t, x))$ of the fluid element at time $t$ and position $x$ and the pressure $p(t, x)$.\\
A translation invariant Banach space of tempered distributions $\mathcal{E}$ is called a critical space for NSE if its norm is invariant under the action of the scaling $f(.) \longrightarrow \lambda f(\lambda .)$. One can take, for example, $\mathcal{E} = L^d(\mathbb{R}^d)$ or the smaller space $\mathcal{E} = \dot{H}^{\frac{d}{2} - 1}(\mathbb{R}^d)$. In fact, one has the “chain of critical spaces” given by the continuous imbedding
\begin{equation}\label{Critical spaces}
\dot{H}^{\frac{d}{2} - 1}(\mathbb{R}^d) \hookrightarrow L^d(\mathbb{R}^d) \hookrightarrow \dot{B}^{\frac{d}{p}-1}_{p, \infty}(\mathbb{R}^d)_{(p < \infty)} \hookrightarrow BMO^{-1}(\mathbb{R}^d) \hookrightarrow \dot{B}^{- 1}_{\infty, \infty}(\mathbb{R}^d).
\end{equation}
It is remarkable feature that the NSE are well-posed in the sense of Hadarmard (existence, uniqueness and continuous dependence on data) when the initial datum is divergence-free and belongs to the critical function spaces (except $\dot{B}^{- 1}_{\infty, \infty}$) listed in $\eqref{Critical spaces}$ (see \cite{M. Cannone 1995} for $\dot{H}^{\frac{d}{2} - 1}(\mathbb{R}^d)$, $L^d(\mathbb{R}^d)$, and $\dot{B}^{\frac{d}{p}-1}_{p, \infty}(\mathbb{R}^d)$,  see \cite{H. Koch 2001} for $BMO^{-1}(\mathbb{R}^d)$, and the recent ill-posedness result  \cite{J. Bourgain 2008} for $\dot{B}^{- 1}_{\infty, \infty}(\mathbb{R}^d)$).\\
In the 1960s, mild solutions were first  constructed by Kato and 
Fujita (\cite{T. Kato 1962}, \cite{H. Fujita 1964}) 
that are continuous in time and take values in the Sobolev space 
$H^s(\mathbb{R}^d),\linebreak  (s \geq \frac{d}{2} - 1)$, say 
$u \in C([0, T]; H^s(\mathbb{R}^d))$. In 1992, a modern treatment 
for mild solutions in $H^s(\mathbb{R}^d), (s \geq \frac{d}{2} - 1)$ 
was given by Chemin \cite{J. M. Chemin 1992}. In 1995, using the simplified 
version of the bilinear operator, Cannone proved the existence 
for mild solutions in $\dot{H}^s(\mathbb{R}^d), (s \geq \frac{d}{2} - 1)$, 
see \cite{M. Cannone 1995}. 
Results on the existence of mild solutions with value in 
$L^p(\mathbb{R}^d), (p > d)$ were established in  the papers 
of  Fabes, Jones and Rivi\`{e}re \cite{E. Fabes 1972a} 
and of Giga \cite{Y. Giga 1986a}. Concerning the initial 
datum in the space $L^\infty$, the existence of a mild solution was 
obtained by Cannone and Meyer in (\cite{M. Cannone 1995}, 
\cite{M. Cannone. Y. Meyer 1995}). Moreover, in (\cite{M. Cannone 1995}, 
\cite{M. Cannone. Y. Meyer 1995}), they also obtained theorems 
on the existence of mild solutions with value in Morrey-Campanato space 
$M^p_2(\mathbb{R}^d), (p > d)$ and Sobolev space $H^s_p(\mathbb{R}^d), (p < d, \frac{1}{p} 
- \frac{s}{d} < \frac{1}{d})$, and in general in the case of a so-called 
well-suited sapce $\mathcal{W}$ for NSE. The NSE in the  
Morrey-Campanato spaces were also treated by Kato  \cite{T. Kato 1992} and 
Taylor \cite{M.E. Taylor 1992}.\\
In 1981, Weissler \cite{F.B. Weissler:1981} gave the first 
existence result of mild solutions in the half space $L^3(\mathbb{R}^3_+)$. 
Then Giga and Miyakawa \cite{Y. Giga 1985} generalized the
result to $L^3(\Omega)$, where $\Omega$ is a bounded 
domain in $\mathbb{R}^3$. Finally, in 1984, Kato \cite{T. Kato 1984} 
obtained, by means of a purely analytical tool (involving only H$\ddot{\text{o}}$lder and 
Young inequalities and without using any estimate of fractional 
powers of the Stokes operator), an existence theorem in the whole 
space $L^3(\mathbb{R}^3)$. In (\cite{M. Cannone 1995}, 
\cite{M. Cannone 1997}, \cite{M. Cannone 1999}), Cannone 
showed how to simplify Kato's proof. The idea is to take advantage 
of the structure of the bilinear operator in its scalar form. 
In particular, the divergence $\nabla$ and heat $e^{t\Delta}$ operators 
can be treated as a single convolution operator. In 1994, Kato and 
Ponce \cite{T. Kato 1994} showed that the NSE are well-posed 
when the initial datum belongs to homogeneous Sobolev spaces 
$\dot{H}^{\frac{d}{p} - 1}_p(\mathbb{R}^d), (d \leq p < \infty)$. 
Recently, the authors of this article have 
considered NSE in mixed-norm Sobolev-Lorentz spaces, 
see \cite{N. M. Tri: Tri2014a}. In  \cite{N. M. Tri: Tri2014???}, 
we showed that NSE are well-posed when  the initial datum belongs 
to Sobolev spaces $\dot{H}^s_p(\mathbb{R}^d)$ with non-positive-regular 
indexes $ (p \geq d, \frac{d}{p}-1 \leq s \leq 0)$. In \cite{N. M. Tri: Tri2014??}, 
we showed that the bilinear operator
\begin{equation}\label{B}
B(u, v)(t) = \int_{0}^{t} e^{(t-\tau ) \Delta} \mathbb{P}
 \nabla.\big(u(\tau)\otimes v(\tau)\big) \dif\tau
\end{equation}
is bicontinuous in $L^\infty([0, T]; \dot{H}^s_p(\mathbb{R}^d))$ 
with super-critical and non-negative-regular indexes $(0 \leq s <d, p > 1,\ {\rm and}\ \frac{s}{d}<\frac{1}{p}<\frac{s+1}{d})$, 
and we established the inequality
\begin{gather*}
\big\|B(u, v)\big\|_{L^\infty([0, T]; \dot{H}^s_p)} \leq
 C_{s,p,d}T^{\frac{1}{2}(1 + s - \frac{d}{p})}
\big\|u\big\|_{L^\infty([0, T]; \dot{H}^s_p)}
\big\|v\big\|_{L^\infty([0, T]; \dot{H}^s_p)}. 
\end{gather*}
In this case existence and uniqueness theorems of local mild 
solutions can therefore be easily deduced. In \cite{N. M. Tri: Tri2014????} 
we prove that NSE are well-posed 
when the initial datum belongs to the Sobolev spaces 
$\dot{H}^{\frac{d}{p} - 1}_p(\mathbb{R}^d)$ with $(1 < p \leq d)$.\\
In this paper, for $s \in \mathbb R$ and $1 \leq p, r \leq \infty$, 
we first recall the notion of the Fourier-Lebesgue spaces $\mathcal{L}^p(\mathbb{R}^d)$,
introduced and investigated in \cite{Lars Hormander 1976}; 
then we introduce and study Sobolev-Fourier-Lebesgue spaces 
$\dot{H}^s_{\mathcal{L}^p}(\mathbb{R}^d)$, 
and  Sobolev-Fourier-Lorentz spaces $\dot{H}^s_{\mathcal{L}^{p, r}}(\mathbb{R}^d)$. 
After that we show that the Navier-Stokes equations 
are well-posed when  the initial datum belongs to the critical 
Sobolev-Fourier-Lorentz spaces $\dot{H}^{\frac{d}{p}-1}_{\mathcal{L}^{p,r}}
(\mathbb{R}^d)$ with $d \geq 2, 1 \leq p < \infty$, and $1 \leq r < \infty$. 
The spaces $\dot{H}^{\frac{d}{p}-1}_{\mathcal{L}^{p,r}}(\mathbb{R}^d)$ 
are more general than the spaces $\dot{H}^{\frac{d}{p}-1}_{\mathcal{L}^p}(\mathbb{R}^d)$. 
In particular,  $\dot{H}^{\frac{d}{p}-1}_{\mathcal{L}^{p,r}}(\mathbb{R}^d) = \dot{H}^{\frac{d}{p}-1}_{\mathcal{L}^p}(\mathbb{R}^d)$ when 
$\frac{1}{p} + \frac{1}{r} = 1$.\\
In 1997, Le Jan and Sznitman \cite{Le Jan  1997} considered 
a very simple space convenient to the study of NSE, 
which is the space E of tempered distributions $f \in \mathcal{S}'(\mathbb{R}^d)$ 
so that $\hat f(\xi)$ is a locally integrable function on  $\mathbb{R}^d$ and
$\underset{\xi}{\rm sup}|\xi|^{d-1}|\hat{f}(\xi)| < \infty$, 
with  $\hat{}$  standing for the Fourier transform. 
This space may be defined as a Besov space based
on the spaces $PM$ of pseudomeasures ($PM$ is the space 
of the image of the Fourier transforms of essentially bounded functions: 
$PM = \mathcal{F}L^{\infty}$). More precisely, 
$E =  \dot{B}^{d-1,\infty}_{PM}(\mathbb{R}^d)$. 
They showed that the bilinear operator $B$ is bicontinuous in 
$L^\infty([0, T]; \dot{B}^{d-1,\infty}_{PM})$ for all $0 < T \leq \infty$. 
Therefore they can easily deduce the existence of global mild 
solutions in spaces $L^\infty([0, \infty); \dot{B}^{d-1,\infty}_{PM})$ 
when norm of the initial value in the spaces $\dot{B}^{d-1,\infty}_{PM}(\mathbb{R}^d)$ 
is small enough. From Definitions \ref{def2.1} and \ref{def2.2} in Section 2, we have 
$$
PM = \mathcal{L}^{1}, \dot{B}^{d-1,\infty}_{PM}(\mathbb{R}^d) 
= \dot{H}^{d-1}_{\mathcal{L}^{1}}(\mathbb{R}^d).
$$

In 2011, Lei and Lin \cite{Zhen Lei Fang 2011} 
showed  that NSE are well-posed when  the initial datum belongs 
to the spaces $\mathcal{X}^{-1}(\mathbb{R}^d)$, which is defined by 
$$
f \in \mathcal{X}^{-1}(\mathbb{R}^d)\ {\rm if\ and\ only\ if}
\ \big\|(-\Delta)^{-\frac{1}{2}}f\big\|_{\mathcal{X}} < \infty, 
{\rm where}\  \big\|f\big\|_{\mathcal{X}} = \big\|\hat f\big\|_{L^1}. 
$$
They established the existence of global mild solutions 
in the space \linebreak  $L^\infty([0, \infty); \mathcal{X}^{-1})$ 
when norm of the initial value in the 
spaces $\mathcal{X}^{-1}(\mathbb{R}^d)$ is small enough.
From Definitions \ref{def2.1} and \ref{def2.2} in Section 2, we see that 
$$
\mathcal{X}^{-1}(\mathbb{R}^d) 
= \dot{H}^{-1}_{\mathcal{L}^{\infty}}(\mathbb{R}^d).
$$
Thus, the spaces $\dot{B}^{d-1,\infty}_{PM}$ and $\mathcal{X}^{-1}$, 
studied in  \cite{Le Jan  1997} and \cite{Zhen Lei Fang 2011}, 
are particular cases of the critical Sobolev-Fourier-Lebesgue spaces 
$\dot{H}^{\frac{d}{p}-1}_{\mathcal{L}^p}$ with $p=1$ and $p=\infty$, respectively.
Note that estimates in the Lorentz spaces were also studied in \cite{C. Ahn:2005a}, 
\cite{H. Hajaiej: 2012a} (see also the refererences therein). 
Very recently, ill-poseness of NSE 
in critical Besov spaces $\dot B^{-1}_{\infty,q}$ 
was investigated in \cite{B. Wang: 2015a}. 

The paper is organized as follows. In Section 2 we introduce and 
investigate the  Sobolev-Fourier-Lorentz spaces and some auxiliary lemmas. 
In Section 3 we present the main results of the paper.   Due to some 
technical difficulties we will consider three cases 
$1 < p \leq d, d \leq q < \infty$, and $p = 1$ separately. In subsection 3.1
we treat the case $1 < p \leq d$. In subsection 3.2 we consider the case 
$d \leq q < \infty$. Finally, in subsection 3.3 we study the case $p = 1$. 
In the sequence, for a space of functions defined on $\mathbb R^d$, 
say $E(\mathbb R^d)$, we will abbreviate it as $E$. Throughout the paper, we sometimes use the notation $A \lesssim B$ as an equivalent to $A \leq CB$ with a uniform constant $C$. The notation $A \simeq B$
means that $A \lesssim B$ and $B \lesssim A$\vskip 1cm 
\centerline{\bf{\S2. SOBOLEV-FOURIER-LORENTZ SPACES}}
\vskip 0.5cm
\begin{Def}\label{def2.1} (Fourier-Lebesgue spaces). 
(See \cite{Lars Hormander 1976}.)\\
For $1 \leq p \leq \infty$, the Fourier-Lebesgue spaces $\mathcal{L}^p(\mathbb{R}^d)$ 
is defined as the space $\mathcal{F}^{-1}(L^{p'}(\mathbb{R}^d)), (\frac{1}{p'} 
+ \frac{1}{p} = 1)$, equipped with the norm 
$$
\big\|f\big\|_{\mathcal{L}^p(\mathbb{R}^d)} :=\big\|\mathcal{F}(f)\big\|_{L^{p'}(\mathbb{R}^d)},
$$
where $\mathcal{F}$ and $\mathcal{F}^{-1}$ denote the Fourier transform and its inverse.
\end{Def} 
\begin{Def}\label{def2.2}(Sobolev-Fourier-Lebesgue spaces).\\
For $s \in \mathbb R$, and $1 \leq p \leq \infty$, the Sobolev-Fourier-Lebesgue spaces $\dot{H}^s_{\mathcal{L}^p}(\mathbb{R}^d)$ is defined as the space  $\dot{\Lambda}^{-s}\mathcal{L}^p(\mathbb{R}^d)$, equipped with the norm 
$$
\big\|u\big\|_{\dot{H}^s_{\mathcal{L}^p}}: = \big\|\dot{\Lambda}^su\big\|_{\mathcal{L}^p}.
$$
where $\dot{\Lambda} = \sqrt{-\Delta}$ is the homogeneous Calderon 
pseudo-differential operator defined as 
$$
\widehat{\dot{\Lambda} g}(\xi) = |\xi|\hat{g}(\xi).
$$ 
\end{Def}
\begin{Def}\label{def2.3} (Lorentz spaces). (See \cite{J. Bergh: 1976}.)\\
For $1 \leq p, r \leq \infty$, the Lorentz space $L^{p, r}(\mathbb{R}^d)$
is defined as follows. \linebreak A measurable function 
$f \in L^{p,r}(\mathbb{R}^d)$ if and only if\\
$\big\|f\big\|_{L^{p,r}}(\mathbb{R}^d) := 
\big(\int_0 ^\infty (t^{\frac{1}{p}}f^*(t))^r
\frac{\dif t}{t}\big)^{\frac{1}{r}} < \infty$ when $1 \leq r < \infty$,\\
$\big\|f\big\|_{L^{p,\infty}}(\mathbb{R}^d) :=  \underset{t > 0}
{\rm sup}\ t^{\frac{1}{p}}f^*(t) < \infty$ when $r = \infty $,\\
 where $f^*(t) = \inf \big\{\tau : \mathcal{M}^d (\{x: |f(x)| > \tau\}) \leq t \big\}$, 
with $\mathcal{M}^d$ being the  Lebesgue measure in $\mathbb R^d$.
\end{Def} 
\begin{Def}\label{def2.4} (Fourier-Lorentz spaces).\\
For $1 \leq p, r \leq \infty$, the Fourier-Lorentz spaces 
$\mathcal{L}^{p, r}(\mathbb{R}^d)$ is defined as the space 
$\mathcal{F}^{-1}(L^{p',r}(\mathbb{R}^d)), (\frac{1}{p'} + \frac{1}{p} = 1)$,
equipped with the norm 
$$
\big\|f\big\|_{\mathcal{L}^{p, r}(\mathbb{R}^d)} 
:=\big\|\mathcal{F}(f)\big\|_{L^{p', r}(\mathbb{R}^d)}.
$$
\end{Def} 
 \begin{Def}\label{def2.5}(Sobolev-Fourier-Lorentz spaces).\\
For $s \in \mathbb R$ and $1 \leq r, p \leq \infty$, the Sobolev-Fourier-Lorentz spaces $\dot{H}^s_{\mathcal{L}^{p, r}}(\mathbb{R}^d)$ is defined as the space  
$\dot{\Lambda}^{-s}\mathcal{L}^{p, r}(\mathbb{R}^d)$, equipped with the norm  
$$
\big\|u\big\|_{\dot{H}^s_{\mathcal{L}^{p, r}}}
: = \big\|\dot{\Lambda}^su\big\|_{\mathcal{L}^{p, r}}.
$$
\end{Def}
\begin{Th}\label{th.2.1} ${\rm (}$Holder's inequality in Fourier-Lorentz spaces${\rm )}$.\\ 
Let $1 <  r, q, \tilde q < \infty$ and $1 \leq h , \tilde h, \hat h \leq +\infty$ satisfy the relations 
$$
\frac{1}{r} = \frac{1}{q} + \frac{1}{\tilde{q}}\ and\ 
 \frac{1}{h} = \frac{1}{\tilde{h}} + \frac{1}{\hat{h}}. 
$$
Suppose that $u \in \mathcal{L}^{q, \tilde{h}}$ and $v \in \mathcal{L}^{\tilde{q}, \hat{h}}$. 
Then $uv \in \mathcal{L}^{r, h}$ and we have the inequality
\begin{equation}\label{th.2.1.1}
\big\|uv\big\|_{\mathcal{L}^{r, h}} \lesssim \big\|u\big\|_{\mathcal{L}^{q, 
\tilde{h}}}\big\|v\big\|_{\mathcal{L}^{\tilde{q}, 
\hat{h}}}.
\end{equation}
\end{Th}
\textbf{Proof}. Let $r', q'$, and $\tilde{q}'$ be such that
$$
\frac{1}{r} + \frac{1}{r'} = 1, \frac{1}{q} + \frac{1}{q'} = 1,
 \ \text{and}\ \frac{1}{\tilde{q}} + \frac{1}{\tilde{q}'} = 1.   
$$
It is easily checked that the following conditions are satisfied
$$
1 < r', q', \tilde{q}' < +\infty\ \text{and}\ \frac{1}{r'} + 1 
= \frac{1}{q'} + \frac{1}{\tilde{q}'}.
$$
We have
\begin{equation}\label{th.2.1.2}
\big\|uv\big\|_{\mathcal{L}^{r, h}} =\big\|\widehat{uv}\big\|_{L^{r', h}} 
= \frac{1}{(2\pi)^{d/2}}\big\|\hat{u}*\hat{v}\big\|_{L^{r', h}}.
\end{equation}
Applying Proposition 2.4 (c) in (\cite{P. G. Lemarie-Rieusset 2002}, p. 20), we have
\begin{equation}\label{th.2.1.3}
\big\|\hat{u}*\hat{v}\big\|_{L^{r', h}} 
\lesssim \big\|\hat{u}\big\|_{L^{q', \tilde{h}}}\big\|\hat{v}\big\|_{L^{\tilde{q}', 
\hat{h}}} = \big\|u\big\|_{\mathcal{L}^{q, \tilde{h}}}\big\|v\big\|_{\mathcal{L}^{\tilde{q}, 
\hat{h}}}.
\end{equation}
Now, the estimate \eqref{th.2.1.1} follows from the equality \eqref{th.2.1.2} 
and the inequality \eqref{th.2.1.3}. \qed
\begin{Th}\label{th.2.2}${\rm (}$Young's inequality for 
convolution in Fourier-Lorentz spaces${\rm )}$. \\
 Let $1 < r, q, \tilde q < \infty$, and $1 \leq h, \tilde h, \hat h \leq \infty$ 
satisfy the relations 
$$
\frac{1}{r} + 1 = \frac{1}{q} + \frac{1}{\tilde{q}}\ and
\ \frac{1}{h} = \frac{1}{\tilde{h}} + \frac{1}{\hat{h}}.
$$
Suppose that $u \in L^{q, \tilde{h}}$ 
and $v \in L^{\tilde{q}, \hat{h}}$. Then $u*v \in L^{r, h}$ 
and the following inequality holds 
\begin{equation}\label{th.2.2.1}
\big\|u*v\big\|_{\mathcal{L}^{r, h}} \lesssim 
\big\|u\big\|_{\mathcal{L}^{q, \tilde{h}}}\big\|v\big\|_{\mathcal{L}^{\tilde{q}, 
\hat{h}}}.
\end{equation}
\end{Th}
\textbf{Proof}.  Let $r', q'$, and $\tilde{q}'$ be such that
$$
\frac{1}{r} + \frac{1}{r'} = 1, \frac{1}{q} + \frac{1}{q'} = 1, 
\ \text{and}\ \frac{1}{\tilde{q}} + \frac{1}{\tilde{q}'} = 1.  
$$
By definition 
\begin{equation}\label{th.2.2.2}
\big\|u*v\big\|_{\mathcal{L}^{r, h}} 
=\big\|\widehat{u*v}\big\|_{L^{r', h}} = (2\pi)^{d/2}\big\|\hat{u}\hat{v}\big\|_{L^{r', h}}. 
\end{equation}
We can check that the following conditions are satisfied
$$
1 < r', q', \tilde{q}' < +\infty\ \text{and}\ \frac{1}{r'} 
= \frac{1}{q'} + \frac{1}{\tilde{q}'}.
$$
Applying Proposition 2.3 (c) in (\cite{P. G. Lemarie-Rieusset 2002}, p. 19), we have
\begin{equation}\label{th.2.2.3}
\big\|\hat{u}\hat{v}\big\|_{L^{r', h}} \lesssim 
\big\|\hat{u}\big\|_{L^{q', \tilde{h}}}\big\|\hat{v}\big\|_{L^{\tilde{q}', 
\hat{h}}} = \big\|u\big\|_{\mathcal{L}^{q, \tilde{h}}}\big\|v\big\|_{\mathcal{L}^{\tilde{q}, 
\hat{h}}}.
\end{equation}
Now, the estimate \eqref{th.2.2.1} follows from the equality 
\eqref{th.2.2.2} and the inequality \eqref{th.2.2.3}. \qed
\begin{Th}\label{th.2.3} {\rm (}Sobolev inequality for 
Sobolev-Fourier-Lorentz spaces{\rm )}.\\
Let $1 < q \leq \tilde{q} <\infty, s, \tilde s \in \mathbb{R}, 
s-\frac{d}{q} = \tilde s -\frac{d}{\tilde q}, \ and \ 1 \leq r \leq \infty $. Then 
\begin{equation}\label {th.2.3.1}
\big\|u\big\|_{\dot{H}^{\tilde s}_{\mathcal{L}^{\tilde q, r}}} 
 \lesssim \big\|u\big\|_{\dot{H}^s_{\mathcal{L}^{q, r}}}, 
\forall u \in \dot{H}^s_{\mathcal{L}^{q, r}}.
\end{equation}
\end{Th}
\textbf{Proof}.  We have
\begin{equation}\label {th.2.3.2}
\big\|u\big\|_{\dot{H}^{\tilde s}_{\mathcal{L}^{\tilde q, r}}} 
= \big\|\dot{\Lambda}^{\tilde s-s}\dot{\Lambda}^su\big\|_{\mathcal{L}^{\tilde{q},r}} 
= \big\| |\xi|^{\tilde s-s}\widehat{\dot{\Lambda}^su}(\xi)\big\|_{L^{\tilde{q}',r}},
\end{equation}
where  
$$
\frac{1}{\tilde{q}} + \frac{1}{\tilde{q}'} = 1.
$$
Note that
$$
|\xi|^{-r} \in L^{\frac{d}{r}, \infty}(\mathbb{R}^d)
\ {\rm for\ all}\ r\ {\rm satisfying} \ 0 < r \leq d.
$$
Applying Proposition 2.3 (c) in (\cite{P. G. Lemarie-Rieusset 2002}, p. 19), we have
\begin{equation}\label{2.10}
\big\| |\xi|^{\tilde s-s}\widehat{\dot{\Lambda}^su}(\xi)\big\|_{L^{\tilde{q}',r}} 
\lesssim \big\| |\xi|^{\tilde s -s}\big\|_{L^{\frac{d}{s-\tilde s}, \infty}}.
\big\|\widehat{\dot{\Lambda}^su}(\xi)\big\|_{L^{q',r}} 
\simeq \big\|u\big\|_{\dot{H}^s_{\mathcal{L}^{q, r}}}.
\end{equation}
The estimate \eqref{th.2.3.1} follows from the equality \eqref{th.2.3.2} 
and the inequality \eqref{2.10}. \qed
\begin{Lem}\label{lem2.1} Let $s \in \mathbb{R},  1 \leq p \leq \infty,
\ and\ 1 \leq r \leq \tilde{r}\leq \infty$.\\
{\rm (a)} We have the following imbedding maps
\begin{gather*}
\mathcal{L}^{p,1} \hookrightarrow \mathcal{L}^{p,r} 
\hookrightarrow \mathcal{L}^{p,\tilde{r}} \hookrightarrow \mathcal{L}^{p,\infty},\\
\dot{H}^s_{\mathcal{L}^{p,1}} \hookrightarrow 
\dot{H}^s_{\mathcal{L}^{p,r}} \hookrightarrow 
\dot{H}^s_{\mathcal{L}^{p,\tilde{r}}} 
\hookrightarrow \dot{H}^s_{\mathcal{L}^{p,\infty}}.
\end{gather*}
{\rm (b)} $\dot{H}^s_{\mathcal{L}^p} = 
\dot{H}^s_{\mathcal{L}^{p,p'}}$ (equality of the norm), 
where $\frac{1}{p} + \frac{1}{p'} = 1$.
\end{Lem}
\textbf{Proof}. It is easily deduced from the properties of the standard Lorentz spaces. \qed
\begin{Lem}\label{lem2.2} Let $s \in \mathbb{R} \ and\ 1 < p < \infty$. We have\\
{\rm (a)} If $1 < q \leq 2$ then $\dot{H}^s_q \hookrightarrow \dot{H}^s_{\mathcal{L}^q}$.\\
{\rm (b)} If $2 \leq q < \infty$ then $\dot{H}^s_{\mathcal{L}^q} \hookrightarrow \dot{H}^s_q$.
\end{Lem}
\textbf{Proof}. It is deduced from Theorem 1.2.1 (\cite{J. Bergh: 1976}, p. 6). \qed
\begin{Lem}\label{lem.2.3} Assume that $1 \leq r, p \leq \infty 
\ and \ k \in \mathbb{N}$, then the two quantities\\
$$
\big\|u\big\|_{\dot{H}^k_{\mathcal{L}^{p, r}}}\ 
and\ \sum_{|\alpha| = k}\big\|\partial^{\alpha}u\big\|_{\mathcal{L}^{p, r}}
$$
are equivalent.
\end{Lem}
\textbf{Proof}. First, we prove that
$$
\sum_{|\alpha| = k}\big\|\partial^{\alpha}u\big\|_{\mathcal{L}^{p, r}} 
\lesssim \big\|u\big\|_{\dot{H}^k_{\mathcal{L}^{p, r}}}.
$$
We have
\begin{gather*}
\sum_{|\alpha| = k}\big\|\partial^{\alpha}u\big\|_{\mathcal{L}^{p, r}} 
= \sum_{|\alpha| = k}\big\| i^k \xi^\alpha \hat u(\xi)\big\|_{L^{p', r}} 
= \sum_{|\alpha| = k}\Big\|\frac{\xi^\alpha}{|\xi|^k}|\xi|^k \hat u(\xi)\Big\|_{L^{p', r}} \\ \leq
 \sum_{|\alpha| = k}\big\||\xi|^k \hat u(\xi)\big\|_{L^{p', r}} 
\lesssim \big\|\widehat{\dot{\Lambda}^ku}(\xi)\big\|_{L^{p',r}} 
= \big\|u\big\|_{\dot{H}^k_{\mathcal{L}^{p, r}}}.
\end{gather*}
Next, we prove that
$$
\big\|u\big\|_{\dot{H}^k_{\mathcal{L}^{p, r}}} 
\lesssim \sum_{|\alpha| = k}\big\|\partial^{\alpha}u\big\|_{\mathcal{L}^{p, r}}.
$$
It is easy to see that for all $\xi \in \mathbb{R}^d$, we have 
$$
|\xi|^k \leq d^{\frac{k}{2}}\sum_{|\alpha| = k}|\xi^\alpha|.
$$
This gives the desired result
\begin{gather*}
\big\|u\big\|_{\dot{H}^k_{\mathcal{L}^{p, r}}} 
= \big\||\xi|^k \hat u(\xi)\big\|_{L^{p', r}} \leq  d^{\frac{k}{2}}
\Big\|\sum_{|\alpha| = k}|\xi^\alpha|\hat u(\xi)\Big\|_{L^{p', r}}\\  
\leq d^{\frac{k}{2}}\sum_{|\alpha| = k}\big\|\xi^\alpha\hat u(\xi)\big\|_{L^{p', r}} 
=  d^{\frac{k}{2}}\sum_{|\alpha| = k}\big\|\partial^{\alpha}u\big\|_{\mathcal{L}^{p, r}}.\qed
\end{gather*}
\begin{Lem}\label{lem.2.4} Let $k \in \mathbb{N}, p \in \mathbb{R}$, 
and $r \in \mathbb{R}$ be such that 
$$
0 \leq k \leq d - 1, \frac{k}{d} < \frac{1}{p} < \frac{1}{2} 
+ \frac{k}{2d},\ and\ 1 \leq r \leq \infty.
$$
Then the following inequality holds
$$
\big\|uv\big\|_{\dot{H}^k_{\mathcal{L}^{q, r}}} 
\lesssim \big\|u\big\|_{\dot{H}^k_{\mathcal{L}^{p, r}}}
\big\|v\big\|_{\dot{H}^k_{\mathcal{L}^{p, r}}}, 
\ \forall  u, v \in \dot{H}^k_{\mathcal{L}^{p, r}},
$$
where 
$$
\frac{1}{q} = \frac{2}{p} - \frac{k}{d}.
$$
\end{Lem}
\textbf{Proof}. First, we estimate 
$\big\|\partial^{\alpha}(uv)\big\|_{\mathcal{L}^{q, r}}$, where
$$
\alpha = (\alpha_1, \alpha_2,...,\alpha_d) \in \mathbb{N}^d, 
\ |\alpha| = \sum_{i = 1}^d\alpha_i = k.
$$
By the general Leibniz rule, we have 
$$
\partial^{\alpha}(uv) = \sum_{\gamma + \beta = \alpha }
\binom{\alpha}{\gamma}(\partial^{\gamma}u)(\partial^{\beta}v).
$$
Set 
$$
\frac{1}{q_1} = \frac{1}{p} - \frac{k - |\gamma|}{d}, \frac{1}{q_2}
 = \frac{1}{p} - \frac{k - |\beta|}{d}.
$$
We have
\begin{gather*}
\frac{1}{q_{1} } + \frac{1}{q_{2}} = \frac{2}{p} - \frac{2k}{d} 
+ \frac{|\gamma| + |\beta|}{d} = \frac{2}{p} - \frac{k}{d} = \frac{1}{q}.
\end{gather*}
Therefore applying Theorems \ref{th.2.1}, \ref{th.2.3}, 
and Lemma \ref{lem2.1} (a) in order to obtain
\begin{gather*}
\big\|(\partial^{\gamma}u)(\partial^{\beta}v)\big\|_{\mathcal{L}^{q, r}} 
\lesssim \big\|\partial^{\gamma}u\big\|_{\mathcal{L}^{q_1, r}}
\big\|\partial^{\beta}v\big\|_{\mathcal{L}^{q_2, \infty}} 
\lesssim \big\|\partial^{\gamma}u\big\|_{\dot{H}^{k - |\gamma|}_{\mathcal{L}^{p, r}}}
\|\partial^{\beta}v\|_{\dot{H}^{k - |\beta|}_{\mathcal{L}^{p, \infty}}}\\
 \lesssim \big\|\partial^{\gamma}u\big\|_{\dot{H}^{k - |\gamma|}_{\mathcal{L}^{p, r}}}
\|\partial^{\beta}v\|_{\dot{H}^{k - |\beta|}_{\mathcal{L}^{p, r}}}
 \lesssim \big\|u\big\|_{\dot{H}^k_{\mathcal{L}^{p, r}}}
\big\|v\big\|_{\dot{H}^k_{\mathcal{L}^{p, r}}}. 
\end{gather*}
Thus, for all $\alpha \in \mathbb{N}^d\ {\rm with} \ |\alpha| = k$, we have 
$$
\big\|\partial^{\alpha}(uv)\big\|_{\mathcal{L}^{q, r}}  
\lesssim \big\|u\big\|_{\dot{H}^k_{\mathcal{L}^{p, r}}}
\big\|v\big\|_{\dot{H}^k_{\mathcal{L}^{p, r}}}.
$$
Applying Lemma \ref{lem.2.3}, we have
$$
\big\|uv\big\|_{\dot{H}^k_{\mathcal{L}^{q, r}}} 
\lesssim \big\|u\big\|_{\dot{H}^k_{\mathcal{L}^{p, r}}}
\big\|v\big\|_{\dot{H}^k_{\mathcal{L}^{p, r}}},
\ \forall  u, v \in \dot{H}^k_{\mathcal{L}^{p, r}}.\qed
$$
\begin{Lem}\label{lem.2.5} Assume that $1 \leq p, r \leq \infty \ and \ s \in \mathbb{R}$. 
If $u_0 \in \dot{H}^s_{\mathcal{L}^{p, r}}$ then $e^{t\Delta}u_0 \in L^\infty([0, \infty); \dot{H}^s_{\mathcal{L}^{p, r}})$ and 
$$
\big\|e^{t\Delta}u_0\big\|_{L^\infty([0, \infty); \dot{H}^s_{\mathcal{L}^{p, r}})} 
\leq\big\|u_0\big\|_{\dot{H}^s_{\mathcal{L}^{p, r}}}.
$$
\end{Lem}
\textbf{Proof}. For $t\ge0$, we have
\begin{gather*}
\big\|e^{t\Delta}u_0\big\|_{\dot{H}^s_{\mathcal{L}^{p, r}}}
 = \big\|e^{t\Delta}\dot{\Lambda}^su_0\big\|_{\mathcal{L}^{p, r}} =
\big\|e^{-t|\xi|^2}|\xi|^s \hat u_0\big\|_{L^{p', r}}  \leq \\
\big\| |\xi|^s \hat u_0\big\|_{L^{p', r}} 
=\big\|\widehat{\dot{\Lambda}^su_0}(\xi)\big\|_{L^{p', r}} 
= \big\|\dot{\Lambda}^su_0(\xi)\big\|_{\mathcal{L}^{p, r}}
=  \big\|u_0\big\|_{\dot{H}^s_{\mathcal{L}^{p, r}}}.\qed 
\end{gather*}
Finally, let us recall the following result on solutions 
of a quadratic equation in Banach spaces 
(Theorem 22.4, \cite{P. G. Lemarie-Rieusset 2002}, p. 227). 
\begin{Th}\label{th.2.4}
Let $E$ be a Banach space, and $B: E \times E \rightarrow  E$ 
be a continuous bilinear form such that there exists $\eta > 0$ so that
$$
\|B(x, y)\| \leq \eta \|x\| \|y\|,
$$
for all x and y in $E$. Then for any fixed $y \in E$ 
such that $\|y\| \leq \frac{1}{4\eta}$, the equation $x = y - B(x,x)$ 
has a unique solution  $\overline{x} \in E$ satisfying $ \|\overline{x}\| \leq \frac{1}{2\eta}$. 
\end{Th}
\vskip 0.5cm
\centerline{\bf{\S3. MAIN RESULTS}} 
\vskip 0.5cm
For $T > 0$, we say that $u$ is a mild solution of NSE on $[0, T]$ 
corresponding to a divergence-free initial data $u_0$ when $u$ 
satisfies the integral equation
$$
u = e^{t\Delta}u_0 - \int_{0}^{t} e^{(t-\tau) \Delta} \mathbb{P} 
\nabla  .\big(u(\tau)\otimes u(\tau)\big) \dif\tau.
$$
Above we have used the following notation: For a tensor $F = (F_{ij})$ 
we define the vector $\nabla.F$ by $(\nabla.F)_i = \sum_{i = 1}^d\partial_jF_{ij}$ 
and for vectors u and v, we define their tensor product $(u \otimes v)_{ij} = u_iv_j$. 
The operator $\mathbb{P}$ is the Leray projection onto the divergence-free fields 
$$
(\mathbb{P}f)_j =  f_j + \sum_{1 \leq k \leq d} R_jR_kf_k, 
$$
where $R_j$ is the Riesz transforms defined as 
$$
R_j = \frac{\partial_j}{\dot \Lambda},\ \ {\rm i.e.} \ \  
\widehat{R_jg}(\xi) = \frac{i\xi_j}{|\xi|}\hat{g}(\xi).
$$
The heat kernel $e^{t\Delta}$ is defined as 
$$
e^{t\Delta}u(x) = ((4\pi t)^{-d/2}e^{-|.|^2/4t}*u)(x).
$$
If $X$ is a normed space and $u = (u_1, u_2,...,u_d), u_i \in X, 1 \leq i \leq d$, then we write 
$$
u \in X, \|u\|_X = \Big(\sum_{i = 1}^d\big\|u_i\big\|_X^2\Big)^{1/2}.
$$
In this main section we investigate mild solutions to NSE when 
the initial datum belongs to critical spaces 
$\dot{H}^{\frac{d}{p} - 1}_{\mathcal{L}^{p,r}}(\mathbb{R}^d)$ 
with $1 \leq p < \infty$ and $1 \leq  r < \infty$. We consider three cases 
$1 < p \leq d, d \leq q < \infty$, and $p = 1$ separately.
\vskip 0.4cm
\textbf{3.1. Solutions to the Navier-Stokes equations with the 
initial value in the critical spaces  
$\dot{H}^{\frac{d}{p} - 1}_{\mathcal{L}^{p,r}}(\mathbb{R}^d)$ 
with $1 < p \leq d$ and $1 \leq r < \infty$}.\\

We define an auxiliary space $\mathcal{K}^{\tilde p}_{p,r,T}$ 
which is made up by the functions $u(t,x)$ such that 
$$
\big\|u\big\|_{\mathcal{K}^{\tilde p}_{p,r,T}}
:= \underset{0 < t < T}{{\rm sup}}t^{\frac{\alpha}{2}}\Big\|u(t,x)
\Big\|_{\dot{H}^{\frac{d}{p} - 1}_{\mathcal{L}^{\tilde p,r}}} < \infty,
$$
and
\begin{equation}\label{3.1.1}
\underset{t \rightarrow 0}{\rm lim\ }t^{\frac{\alpha}{2}}
\Big\|u(t,x)\Big\|_{\dot{H}^{\frac{d}{p} - 1}_{\mathcal{L}^{\tilde p,r}}} = 0,
\end{equation}
with 
$$
1 < p \leq \tilde p < \infty, \frac{1}{p}-\frac{1}{d} < \frac{1}{\tilde p}, 1 \leq r \leq \infty, T >0,
$$
and 
$$
\alpha = \alpha(p,\tilde p) = d\Big(\frac{1}{p} - \frac{1}{\tilde p}\Big).
$$
In the case $\tilde p = p$, it is also convenient to 
define the space $\mathcal{K}^p_{p,r,T}$ as the natural 
space  $L^\infty\big([0, T]; \dot{H}^{\frac{d}{p} - 1}_{\mathcal{L}^{p,r}}\big)$ 
with the additional condition that its elements $u(t,x)$ satisfy 
\begin{equation}\label{3.1.2}
\underset{t \rightarrow 0}{\text{lim}\ }
\Big\|u(t,x)\Big\|_{\dot{H}^{\frac{d}{p} - 1}_{\mathcal{L}^{p,r}}} = 0.
\end{equation}
\begin{Lem}\label{lem3.1} Let $1 \leq r \leq \tilde{r}\leq \infty$. 
Then we have the following imbedding
$$
\mathcal{K}^{\tilde p}_{p,1,T} \hookrightarrow \mathcal{K}^{\tilde p}_{p,r,T}
 \hookrightarrow \mathcal{K}^{\tilde p}_{p,\tilde r,T} \hookrightarrow 
\mathcal{K}^{\tilde p}_{p,\infty,T}.
$$
\end{Lem}
\textbf{Proof}. It is easily deduced from Lemma \ref{lem2.1} (a) 
and the definition \linebreak of  $\mathcal{K}^{\tilde p}_{p,r,T}$. \qed

\begin{Lem}\label{lem.3.2} Suppose that $u_0 \in \dot{H}^{\frac{d}{p} - 1}_{\mathcal{L}^{p,r}}(\mathbb{R}^d)$ with $1 < p \leq d$ and $1 \leq r < \infty$, 
then $e^{t\Delta}u_0 \in \mathcal{K}^{\tilde p}_{p,1,\infty}$ with $\frac{1}{p} - \frac{1}{d} < \frac{1}{\tilde p} < \frac{1}{p}$.
\end{Lem}
\textbf{Proof}. Before proving this lemma, we need to prove the following lemma.
\begin{Lem}\label{lem.3.3} Suppose that 
$u_0 \in  L^{q, r}(\mathbb{R}^d)$ with $1 \leq q \leq \infty$ and $1 \leq r < \infty$. Then 
$\underset{n \rightarrow \infty}{\rm lim}\big\|1_{B^c_n}u_0\big\|_{L^{q, r}} = 0$,
where $n \in \mathbb{N}, B_n = \{x \in \mathbb{R}^d : |x| < n\}, B^c_n = \mathbb{R}^d \backslash B_n$, and $1_{B^c_n}$ 
is the indicator function of the set $B^c_n$ on $\mathbb{R}^d:  1_{B^c_n}(x) = 1$ for $x \in B^c_n$ and $1_{B^c_n}(x) = 0$ otherwise.
\end{Lem}
\textbf{Proof}.  With $\delta > 0$ being fixed, we have
\begin{equation}\label{lem.3.3.1}
\big\{x: |1_{B^c_n}u_0(x)| > \delta\big\} \supseteq \big\{x: |1_{B^c_{n + 1}}u_0(x)| > \delta\big\},
\end{equation}
and
\begin{equation}\label{lem.3.3.2}
\underset{n = 0}{\overset{\infty }{\cap}}\{x: |1_{B^c_n}u_0(x)| > \delta\} = \emptyset.
\end{equation}
Note that 
$$
\mathcal{M}^d\big(\{x: |1_{B^c_0}u_0(x)| > \delta\}\big) = \mathcal{M}^d\big(\{x: |u_0(x)| > \delta\}\big).
$$
We prove that
\begin{equation}\label{lem.3.3.3}
\mathcal{M}^d\big(\{x: |u_0(x)| > \delta\}\big) < \infty,
\end{equation}
assuming on the contrary
$$
\mathcal{M}^d\big(\{x: |u_0(x)| > \delta\}\big) = \infty.
$$
Set
$$
u_0^*(t) = \inf \big\{\tau : \mathcal{M}^d\big(\{x: |u_0(x)| > \tau\}\big) \leq t \big\}.
$$
We have $u_0^*(t) \geq \delta\ \text{for all}\ t > 0$, from the definition 
of the Lorentz space, we get 
$$
\big\|u_0\big\|_{L^{q,r}} = \Big(\int_0 ^\infty (t^{\frac{1}{q}}u_0^*(t))^r
\frac{\dif t}{t}\Big)^{\frac{1}{r}} \geq \Big(\int_0 ^\infty (t^{\frac{1}{q}}\delta)^r
\frac{\dif t}{t}\Big)^{\frac{1}{r}} =\delta \Big(\int_0 ^\infty t^{\frac{r}{q} - 1}
\dif t\Big)^{\frac{1}{r}} = \infty,
$$
a contradiction.\\
From \eqref{lem.3.3.1}, \eqref{lem.3.3.2}, and \eqref{lem.3.3.3}, we have
\begin{equation}\label{lem.3.3.4}
\underset{n \rightarrow \infty}{\rm lim}\mathcal{M}^d
\big(\{x: |1_{B^c_n}u_0(x)| > \delta\}\big) = 0.
\end{equation}
Set
$$
u^*_n(t) = \inf \big\{\tau : \mathcal{M}^d
\big(\{x: |1_{B^c_n}u_0(x)| > \tau\}\big) \leq t \big\}.
$$
We have
\begin{equation}\label{lem.3.3.5}
u^*_n(t) \geq u^*_{n + 1}(t).
\end{equation}
Fixed $t > 0$. For any $\epsilon > 0$, from \eqref{lem.3.3.4} it follows that 
there exist \linebreak $n_0 = n_0(t, \epsilon)$ is large enough such that 
$$
\mathcal{M}^d\big(\{x: |1_{B^c_n}u_0(x)| 
> \epsilon\}\big) \leq t , \forall n \geq n_0.
$$
From this we deduce that
$$
u^*_n(t) \leq \epsilon , \forall n \geq n_0,
$$
therefore
\begin{equation}\label{lem.3.3.6}
\underset{n \rightarrow \infty}{\text{lim}u^*_n(t) = 0}.
\end{equation}
From \eqref{lem.3.3.5} and \eqref{lem.3.3.6}, we apply 
Lebesgue's monotone convergence theorem to get
$$
\underset{n \rightarrow \infty}{\rm lim}\big\|1_{B^c_n}u_0\big\|_{L^{q, r}} =
 \underset{n \rightarrow \infty}{\rm lim}\Big(\int_0 ^\infty 
(t^{\frac{1}{q}}u^*_n(t))^r\frac{\dif t}{t}\Big)^{\frac{1}{r}} = 0. \qed
$$

Now we return to prove Lemma \ref{lem.3.2}. We prove that
\begin{equation}\label{lem.3.3.7}
\underset{0 < t < \infty}{\rm sup}t^{\frac{\alpha}{2}}
\Big\|e^{t\Delta}u_0\Big\|_{\dot{H}^{\frac{d}{p} - 1}_{\mathcal{L}^{\tilde p,1}}}
 \lesssim \Big\|u_0\Big\|_{\dot{H}^{\frac{d}{p} - 1}_{\mathcal{L}^{p,r}}}.
\end{equation}
Let $p'$ and $\tilde p'$ be such that
$$
\frac{1}{p} + \frac{1}{p'} = 1\ {\rm and}\ \frac{1}{\tilde p} + \frac{1}{\tilde p'} = 1.
$$
We have
\begin{equation}\label{lem.3.3.8}
\Big\|e^{t\Delta}u_0\Big\|_{\dot{H}^{\frac{d}{p} - 1}_{\mathcal{L}^{\tilde p,1}}} =
\big\|e^{-t\left| \xi \right|^2}|\xi|^{\frac{d}{p} - 1}\hat u_0(\xi)\big\|_{L^{\tilde p', 1}_\xi}.
\end{equation}
Applying Holder's inequality in the Lorentz spaces (see Proposition 2.3 (c) 
in \cite{P. G. Lemarie-Rieusset 2002}, p. 19), we have
\begin{gather}
\big\|e^{-t\left| \xi \right|^2}|\xi|^{\frac{d}{p} - 1}\hat u_0(\xi)\big\|_{L^{\tilde p', 1}_\xi} 
\lesssim \big\|e^{-t\left| \xi \right|^2}\big\|_{L^{\frac{p\tilde p}{\tilde p-p}, 1}_{\xi}}
\big\||\xi|^{\frac{d}{p} - 1}\hat u_0(\xi)\big\|_{L^{p', \infty}}= \notag\\
t^{-\frac{d}{2}(\frac{1}{p} - \frac{1}{\tilde p})}\big\|e^{-\left| \xi 
\right|^2}\big\|_{L^{\frac{p\tilde p}{\tilde p-p}, 1}}
\big\||\xi|^{\frac{d}{p} - 1}\hat u_0(\xi)\big\|_{L^{p', \infty}} \lesssim 
 t^{-\frac{\alpha}{2}}\big\||\xi|^{\frac{d}{p} - 1}\hat u_0(\xi)\big\|_{L^{p', r}}\notag \\ 
= t^{-\frac{\alpha}{2}}\Big\|u_0\Big\|_{\dot{H}^{\frac{d}{p} - 1}_{\mathcal{L}^{p,r}}}. \label{lem.3.3.9}
\end{gather}
The estimate \eqref{lem.3.3.7} follows from the 
equality \eqref{lem.3.3.8} and the estimate \eqref{lem.3.3.9}. \\
We claim now that
\begin{equation}\label{lem.3.3.10}
\underset{t \rightarrow 0}{\rm lim\ }t^{\frac{\alpha}{2}}
\Big\|e^{t\Delta}u_0\Big\|_{\dot{H}^{\frac{d}{p} - 1}_{\mathcal{L}^{\tilde p,1}}} = 0.
\end{equation}
From the equality \eqref{lem.3.3.8}, we have
\begin{gather*}
t^{\frac{\alpha}{2}}\Big\|e^{t\Delta}u_0\Big\|_{\dot{H}^{\frac{d}{p} 
- 1}_{\mathcal{L}^{\tilde p,1}}} \leq 
t^{\frac{\alpha}{2}}\big\|e^{-t\left| \xi \right|^2}1_{B_n^c}|\xi|^{\frac{d}{p}-1}
\hat u_0(\xi)\big\|_{L^{\tilde p', 1}_\xi} \\
+\ t^{\frac{\alpha}{2}}\big\|e^{-t\left| \xi \right|^2}1_{B_n}|\xi|^{\frac{d}{p}-1}
\hat u_0(\xi)\big\|_{L^{\tilde p', 1}_\xi}.
\end{gather*}
For any $\epsilon > 0$. Applying Holder's inequality in the Lorentz spaces and 
using Lemma \ref{lem.3.3}, we have
\begin{gather}
t^{\frac{\alpha}{2}}\big\|e^{-t\left| \xi \right|^2}1_{B^c_n}|\xi|^{\frac{d}{p}-1}
\hat u_0(\xi)\big\|_{L^{\tilde p', 1}_\xi} \leq 
Ct^{\frac{\alpha}{2}}\big\|e^{-t\left| \xi \right|^2}\big\|_{L^{\frac{p\tilde p}{\tilde p-p}, 1}_{\xi}}
\big\|1_{B^c_n}|\xi|^{\frac{d}{p}-1}\hat u_0(\xi)\big\|_{L^{p', \infty}}= \notag \\
C\big\|e^{-\left| \xi \right|^2}\big\|_{L^{\frac{p\tilde p}{\tilde p-p}, 1}}
\big\|1_{B^c_n}|\xi|^{\frac{d}{p}-1}\hat u_0(\xi)\big\|_{L^{p', \infty}} 
\leq C'\big\|1_{B^c_n}|\xi|^{\frac{d}{p}-1}\hat u_0(\xi)\big\|_{L^{p', r}} 
< \frac{\epsilon}{2} \ \  \label{lem.3.3.11}
\end{gather}
for large enough $n$. Fixed one of such $n$ and applying Holder's inequality in the 
Lorentz spaces, we have
\begin{gather}
t^{\frac{\alpha}{2}}\big\|e^{-t\left| \xi \right|^2}1_{B_n}|\xi|^{\frac{d}{p}-1}
\hat u_0(\xi)\big\|_{L^{\tilde p', 1}_\xi} \leq 
Ct^{\frac{\alpha}{2}}\big\|1_{B_n}e^{-t\left|\xi \right|^2}\big\|_{L^{\frac{p\tilde p}{\tilde p - p}, 1}_{\xi}}
\big\||\xi|^{\frac{d}{p}-1}\hat u_0(\xi)\big\|_{L^{p', \infty}} \notag \\
\leq  Ct^{\frac{\alpha}{2}}\big\|1_{B_n}\big\|_{L^{\frac{p\tilde p}{\tilde p - p}, 1}}
\big\||\xi|^{\frac{d}{p}-1}\hat u_0(\xi)\big\|_{L^{p', \infty}}
\leq C''(n)t^{\frac{\alpha}{2}}\big\||\xi|^{\frac{d}{p}-1}\hat u_0(\xi)\big\|_{L^{p', r}} \notag \\
=  C''(n)t^{\frac{\alpha}{2}}\Big\|u_0\Big\|_{\dot{H}^{\frac{d}{p}
-1}_{\mathcal{L}^{p, r}}} < \frac{\epsilon}{2} \label{lem.3.3.12}
\end{gather}
for small enough $t =t(n) > 0$. From estimates \eqref{lem.3.3.11} 
and \eqref{lem.3.3.12}, we have,
$$
t^{\frac{\alpha}{2}}\Big\|e^{t\Delta}u_0\Big\|_{\dot{H}^{\frac{d}{p} 
- 1}_{\mathcal{L}^{\tilde p,1}}} \leq C'\big\|1_{B^c_n}|\xi|^{\frac{d}{p}-1}
\hat u_0(\xi)\big\|_{L^{p', r}} + 
C''(n)t^{\frac{\alpha}{2}}\Big\|u_0\Big\|_{\dot{H}^{\frac{d}{p}
-1}_{\mathcal{L}^{p, r}}} < \epsilon. \qed
$$
In the following lemmas a particular attention will be devoted to the study 
of the bilinear operator $B(u, v)(t)$ defined by \eqref{B}.\\
In the following lemmas, denote by $[x]$ the integer part of $x$ 
and by $\{x\}$ the fraction part of $x$. 
\begin{Lem}\label{lem.3.4}
Let $1 <p\leq d$. Then for all $\tilde p$ be such that
\begin{equation}\label{lem.3.4.1}
\frac{1}{2p}+\frac{[\frac{d}{p}]-1}{2d} < \frac{1}{\tilde p} < {\rm min}\Big\{\frac{[\frac{d}{p}]}{d},\frac{1}{2}
+\frac{[\frac{d}{p}]-1}{2d}\Big\},
\end{equation}
the bilinear operator $B(u, v)(t)$ is continuous from 
$\mathcal{K}^{\tilde p}_{\frac{d}{[\frac{d}{p}]},\infty,T} \times \mathcal{K}^{\tilde p}_{\frac{d}{[\frac{d}{p}]},\infty,T}$ 
into $\mathcal{K}^p_{p,1,T}$ and the following inequality holds 
\begin{equation}\label{lem.3.4.2}
\big\|B(u, v)\big\|_{\mathcal{K}^p_{p,1,T}} \leq
 C\Big\|u\Big\|_{\mathcal{K}^{\tilde p}_{\frac{d}{[\frac{d}{p}]},\infty,T}}
\Big\|v\Big\|_{\mathcal{K}^{\tilde p}_{\frac{d}{[\frac{d}{p}]},\infty,T}},
\end{equation}
where C is a positive constant and independent of T.
\end{Lem}
\textbf{Proof}. We have
\begin{gather}
\Big\|B(u, v)(t)\Big\|_{\dot{H}^{\frac{d}{p}-1}_{\mathcal{L}^{p, 1}}}
 \leq \int_{0}^{t} \Big\|e^{(t-\tau ) \Delta} \mathbb{P} 
\nabla .(u(\tau)\otimes v(\tau))\Big\|_{\dot{H}^{\frac{d}{p}-1}_{\mathcal{L}^{p, 1}}} \dif\tau\notag \\
= \int_{0}^{t} \Big\|\dot{\Lambda}^{\frac{d}{p}-1}e^{(t-\tau ) \Delta} \mathbb{P} \nabla .(u(\tau)
\otimes v(\tau))\Big\|_{\mathcal{L}^{p, 1}}\dif\tau.\label{lem.3.4.3}
\end{gather}
Note that 
\begin{gather*}
\Big(\dot{\Lambda}^{\frac{d}{p}-1}e^{(t - \tau)\Delta}\mathbb{P}\nabla .\big(u(\tau) 
\otimes v(\tau)\big)\Big)_j^{\wedge}(\xi)  \\ 
=\Big(\dot{\Lambda}^{\{\frac{d}{p}\}}e^{(t - \tau)\Delta}\mathbb{P}
\nabla .\dot{\Lambda}^{[\frac{d}{p}]-1}\big(u(\tau) 
\otimes v(\tau)\big)\Big)_j^{\wedge}(\xi) \\
= |\xi|^{\{\frac{d}{p}\}}e^{-(t - \tau)|\xi|^2}\sum_{l, k = 1}^d \Big(\delta_{jk} - 
\frac{\xi_j\xi_k}{|\xi|^2}\Big)(i\xi_l)\Big(
\dot{\Lambda}^{[\frac{d}{p}]-1}\big(u_l (\tau)v_k(\tau)\big)\Big)^{\wedge}(\xi).
\end{gather*}
Thus
\begin{gather}
\Big(\dot{\Lambda}^{\frac{d}{p}-1}e^{(t - \tau)\Delta}\mathbb{P}
\nabla .\big(u(\tau) \otimes v(\tau)\big)\Big)_j  \notag \\
=\frac{1}{(t -\tau)^{\frac{\{\frac{d}{p}\} + d+1}{2}}}\sum_{l, k = 1}^d K_{l, k, j}
\Big(\frac{.}{\sqrt{t - \tau}}\Big)*\Big(\dot{\Lambda}^{[\frac{d}{p}]-1}
\big(u_l (\tau)v_k(\tau)\big)\Big),\label{lem.3.4.4}
\end{gather}
where 
\begin{gather}
\widehat{K_{l, k, j}}(\xi) = \frac{1}{(2\pi)^{d/2}}|\xi|^{\{\frac{d}{p}\}}e^{-|\xi|^2}
\Big(\delta_{jk} - \frac{\xi_j\xi_k}{|\xi|^2}\Big)(i\xi_l).\label{lem.3.4.5}
\end{gather}
Setting the tensor $K(x) = \{K_{l, k, j}(x)\}$, we can rewrite 
the equality \eqref{lem.3.4.4} in the tensor form 
\begin{gather*}
\dot{\Lambda}^{\frac{d}{p}-1}e^{(t - \tau)\Delta}\mathbb{P}
\nabla .\big(u(\tau) \otimes v(\tau)\big) \notag \\
=\frac{1}{(t - \tau)^{\frac{\{\frac{d}{p}\}+d+1}{2}}}
K\Big(\frac{.}{\sqrt{t - \tau}}\Big)*\Big(\dot{\Lambda}^{[\frac{d}{p}]
-1}\big(u(\tau) \otimes v(\tau)\big)\Big).
\end{gather*}
Applying Theorem \ref{th.2.2} for convolution in the Fourier-Lorentz spaces, we have
\begin{gather}
\Big\|\dot{\Lambda}^{\frac{d}{p}-1}e^{(t - \tau)\Delta}\mathbb{P}
\nabla .\big(u(\tau) \otimes v(\tau)\big)\Big\|_{\mathcal{L}^{p, 1}} \lesssim \notag \\
\frac{1}{(t - \tau)^{\frac{\{\frac{d}{p}\}+d+1}{2}}}\Big\|
K\Big(\frac{.}{\sqrt{t - \tau}}\Big)\Big\|_{\mathcal{L}^{r, 1}}
\big\|\dot{\Lambda}^{[\frac{d}{p}]-1}\big(u(\tau) \otimes v(\tau)\big)
\big\|_{\mathcal{L}^{q, \infty}},\label{lem.3.4.7}
\end{gather}
where
\begin{gather}
\frac{1}{q} = \frac{2}{\tilde p}-\frac{[\frac{d}{p}]-1}{d}\ {\rm and}
\  \frac{1}{r} = 1+ \frac{1}{p} - \frac{2}{\tilde p}+\frac{[\frac{d}{p}]-1}{d}.\label{lem.3.4.8}
\end{gather}
Note that from the inequality \eqref{lem.3.4.1}, we can check 
that $r$ and $q$ satisfy the relations  
$$
1 < r, q < \infty, \frac{1}{p}+1 = \frac{1}{q}+\frac{1}{r}.
$$
Applying Lemma \ref{lem.2.4}, we have
\begin{gather}
\big\| \dot{\Lambda}^{[\frac{d}{p}]-1}\big(u(\tau) \otimes v(\tau)\big)
\big\|_{\mathcal{L}^{q, \infty}} \lesssim \big\|u(\tau)
\big\|_{\dot{H}^{[\frac{d}{p}]-1}_{\mathcal{L}^{\tilde p, \infty}}}
\big\|v(\tau)\big\|_{\dot{H}^{[\frac{d}{p}]-1}_{\mathcal{L}^{\tilde p, \infty}}}.
\label{lem.3.4.9}
\end{gather}
From the equalities \eqref{lem.3.4.5} and \eqref{lem.3.4.8}, we obtain 
\begin{gather}
\Big\|K\Big(\frac{.}{\sqrt{t - \tau}}\Big)\Big\|_{\mathcal{L}^{r,1}} 
= (t - \tau)^{\frac{d}{2}}\Big\|\hat K\big(\sqrt{t - \tau}\ (.)\big)\Big\|_{L^{r',1}}=\notag \\
(t - \tau)^{\frac{d}{2}- \frac{d}{2.r'}}\big\|\hat K\big\|_{L^{r', 1}} 
= (t - \tau)^{\frac{d}{2.r}}\big\|\hat K\big\|_{L^{r', 1}} 
\simeq (t - \tau)^{\frac{d}{2}\big(1+ \frac{1}{p} - \frac{2}{\tilde p}+\frac{[\frac{d}{p}]-1}{d}\big)}.
\label{lem.3.4.10}
\end{gather}
From the estimates \eqref{lem.3.4.7}, \eqref{lem.3.4.9}, 
and \eqref{lem.3.4.10}, we deduce that 
\begin{gather*}
\Big\|\dot{\Lambda}^{\frac{d}{p}-1}e^{(t - \tau)\Delta}\mathbb{P}
\nabla .\big(u(\tau) \otimes v(\tau)\big)\Big\|_{\mathcal{L}^{p, 1}}
 \lesssim  (t - \tau)^{[\frac{d}{p}]-\frac{d}{\tilde p}-1}
\Big\|u(\tau)\Big\|_{\dot{H}^{[\frac{d}{p}]-1}_{\mathcal{L}^{\tilde p, \infty}}}
\Big\|v(\tau)\Big\|_{\dot{H}^{[\frac{d}{p}]-1}_{\mathcal{L}^{\tilde p, \infty}}}\\
=(t - \tau)^{\alpha-1}\Big\|u(\tau)\Big\|_{\dot{H}^{[\frac{d}{p}]
-1}_{\mathcal{L}^{\tilde p, \infty}}}\Big\|v(\tau)\Big\|_{\dot{H}^{[\frac{d}{p}]
-1}_{\mathcal{L}^{\tilde p, \infty}}},
\end{gather*}
where
$$
\alpha = \alpha\Big (\frac{d}{[\frac{d}{p}]},\tilde p\Big) 
= \Big[\frac{d}{p}\Big]-\frac{d}{\tilde p}\ ,
$$
this gives the desired result
\begin{gather*}
\Big\|B(u, v)(t)\Big\|_{\dot{H}^{\frac{d}{p}-1}_{\mathcal{L}^{p, 1}}} 
\lesssim  \int_0^t (t - \tau)^{\alpha-1}\Big\|u(\tau)\Big\|_{\dot{H}^{[\frac{d}{p}]
-1}_{\mathcal{L}^{\tilde p, \infty}}}\Big\|v(\tau)\Big\|_{\dot{H}^{[\frac{d}{p}]
-1}_{\mathcal{L}^{\tilde p, \infty}}}\dif\tau \\
\lesssim \int_0^t (t - \tau)^{\alpha - 1}\tau^{-\alpha}\underset{0 < \eta < t}
{\rm sup}\eta^{\frac{\alpha}{2}}\Big\|u(\eta)
\Big\|_{\dot{H}^{[\frac{d}{p}]-1}_{\mathcal{L}^{\tilde p, \infty}}}
\underset{0 < \eta < t}{\rm sup}
\eta^{\frac{\alpha}{2}}\Big\|v(\eta)\Big\|_{\dot{H}^{[\frac{d}{p}]
-1}_{\mathcal{L}^{\tilde p, \infty}}}\dif\tau 
\end{gather*}
\begin{gather}
= \underset{0 < \eta < t} {\rm sup}\eta^{\frac{\alpha}{2}}\Big\|u(\eta)
\Big\|_{\dot{H}^{[\frac{d}{p}]-1}_{\mathcal{L}^{\tilde p, \infty}}}
\underset{0 < \eta < t}{\rm sup}
\eta^{\frac{\alpha}{2}}\Big\|v(\eta)\Big\|_{\dot{H}^{[\frac{d}{p}]
-1}_{\mathcal{L}^{\tilde p, \infty}}}\int_0^t (t - \tau)^{\alpha - 1}\tau^{-\alpha}\dif\tau \notag \\
\simeq \underset{0 < \eta < t}{\rm sup}\eta^{\frac{\alpha}{2}}
\Big\|u(\eta)\Big\|_{\dot{H}^{[\frac{d}{p}]-1}_{\mathcal{L}^{\tilde p, \infty}}}
\underset{0 < \eta < t}{\rm sup}\eta^{\frac{\alpha}{2}}
\Big\|v(\eta)\Big\|_{\dot{H}^{[\frac{d}{p}]-1}_{\mathcal{L}^{\tilde p, \infty}}}. \label{lem.3.4.11}
\end{gather}
Let us now check the validity of the condition \eqref{3.1.2} for the bilinear term $B(u,v)(t)$. 
Indeed, from \eqref{lem.3.4.11}
$$
\underset{t \rightarrow 0}{\rm lim}\Big\|B(u, v)(t)
\Big\|_{\dot{H}^{\frac{d}{p}-1}_{\mathcal{L}^{p, 1}}}  = 0,
$$
whenever
$$
\underset{t \rightarrow 0}{\rm lim\ }t^{\frac{\alpha}{2}}
\Big\|u(t)\Big\|_{\dot{H}^{[\frac{d}{p}]-1}_{\mathcal{L}^{\tilde p, \infty}}} 
= \underset{t \rightarrow 0}{\rm lim\ }t^{\frac{\alpha}{2}}
\Big\|v(t)\Big\|_{\dot{H}^{[\frac{d}{p}]-1}_{\mathcal{L}^{\tilde p, \infty}}} = 0.
$$
The estimate  \eqref{lem.3.4.2} is deduced from the inequality \eqref{lem.3.4.11}. \qed
\begin{Lem}\label{lem.3.5}
Let $1 <p\leq d$. Then for all $\tilde p$ be such that
\begin{equation}\label{lem.3.5.1}
\frac{[\frac{d}{p}]-1}{d} < \frac{1}{\tilde p} < {\rm min}
\Big\{\frac{[\frac{d}{p}]}{d},\frac{1}{2}+\frac{[\frac{d}{p}]-1}{2d}\Big\},
\end{equation}
the bilinear operator $B(u, v)(t)$ is continuous from 
$\mathcal{K}^{\tilde p}_{\frac{d}{[\frac{d}{p}]},\infty,T} \times \mathcal{K}^{\tilde p}_{\frac{d}{[\frac{d}{p}]},\infty,T}$ 
into $\mathcal{K}^{\tilde p}_{\frac{d}{[\frac{d}{p}]},1,T}$ 
and the following inequality holds
\begin{equation}\label{lem.3.5.2}
\big\|B(u, v)\Big\|_{\mathcal{K}^{\tilde p}_{\frac{d}{[\frac{d}{p}]},1,T}} \leq
 C\Big\|u\Big\|_{\mathcal{K}^{\tilde p}_{\frac{d}{[\frac{d}{p}]},\infty,T}}
\Big\|v\Big\|_{\mathcal{K}^{\tilde p}_{\frac{d}{[\frac{d}{p}]},\infty,T}},
\end{equation}
where C is a positive constant and independent of T.
\end{Lem}
\textbf{Proof}. First, arguing as in Lemma \ref{lem.3.4}, we derive
\begin{gather*}
\dot{\Lambda}^{[\frac{d}{p}]-1}e^{(t - \tau)\Delta}\mathbb{P}
\nabla .\big(u(\tau) \otimes v(\tau)\big) \notag \\
=\frac{1}{(t - \tau)^{\frac{d+1}{2}}}K\Big(\frac{.}{\sqrt{t - \tau}}
\Big)*\Big(\dot{\Lambda}^{[\frac{d}{p}]-1}\big(u(\tau) \otimes v(\tau)\big)\Big),
\end{gather*}
where
\begin{gather}
\widehat{K_{l, k, j}}(\xi) = \frac{1}{(2\pi)^{d/2}}e^{-|\xi|^2}
\Big(\delta_{jk} - \frac{\xi_j\xi_k}{|\xi|^2}\Big)(i\xi_l).\label{eq1}
\end{gather}
Applying Theorem \ref{th.2.2} for the convolution in 
the Fourier-Lorentz spaces, we have
\begin{gather}
\Big\|\dot{\Lambda}^{[\frac{d}{p}]-1}e^{(t - \tau)\Delta}\mathbb{P}
\nabla .\big(u(\tau) \otimes v(\tau)\big)\Big\|_{\mathcal{L}^{\tilde p, 1}}\notag \\
 \lesssim \frac{1}{(t - \tau)^{\frac{d+1}{2}}}
\Big\|K\Big(\frac{.}{\sqrt{t - \tau}}\Big)\Big\|_{\mathcal{L}^{r, 1}}
\big\|\dot{\Lambda}^{[\frac{d}{p}]-1}\big(u(\tau) 
\otimes v(\tau)\big)\big\|_{\mathcal{L}^{q, \infty}},\label{lem.3.5.6}
\end{gather}
where
\begin{gather}
\frac{1}{q} = \frac{2}{\tilde p}-\frac{[\frac{d}{p}]-1}{d}\ {\rm and}
\  \frac{1}{r} = 1- \frac{1}{\tilde p}+\frac{[\frac{d}{p}]-1}{d}.\label{lem.3.5.7}
\end{gather}
Note that from the inequality \eqref{lem.3.5.1}, we can check 
that $r$ and $q$ satisfy the relations  
$$
1 < r, q < \infty, \frac{1}{p}+1 = \frac{1}{q}+\frac{1}{r}.
$$
Applying Lemma \ref{lem.2.4}, we have
\begin{gather}
\big\| \dot{\Lambda}^{[\frac{d}{p}]-1}\big(u(\tau) 
\otimes v(\tau)\big)\big\|_{\mathcal{L}^{q, \infty}} 
\lesssim \Big\|u(\tau)\Big\|_{\dot{H}^{[\frac{d}{p}]
-1}_{\mathcal{L}^{\tilde p, \infty}}}\Big\|v(\tau)
\Big\|_{\dot{H}^{[\frac{d}{p}]-1}_{\mathcal{L}^{\tilde p, \infty}}}.
\label{lem.3.5.8}
\end{gather}
From the equalities \eqref{eq1} and \eqref{lem.3.5.7}, we obtain 
\begin{gather}
\Big\|K\Big(\frac{.}{\sqrt{t - \tau}}\Big)\Big\|_{\mathcal{L}^{r,1}} 
=  (t - \tau)^{\frac{d}{2.r}}\big\|\hat K\big\|_{L^{r', 1}} 
\simeq (t - \tau)^{\frac{d}{2}\big(1 - \frac{1}{\tilde p}+\frac{[\frac{d}{p}]-1}{d}\big)}.
\label{lem.3.5.9}
\end{gather}
From the estimates \eqref{lem.3.5.6}, \eqref{lem.3.5.8}, 
and \eqref{lem.3.5.9}, we deduce that 
\begin{gather*}
\Big\|\dot{\Lambda}^{[\frac{d}{p}]-1}e^{(t - \tau)\Delta}\mathbb{P}
\nabla .\big(u(\tau) \otimes v(\tau)\big)\Big\|_{\mathcal{L}^{\tilde p, 1}}  \\ 
\lesssim  (t - \tau)^{\frac{1}{2}([\frac{d}{p}]-\frac{d}{\tilde p})-1}
\Big\|u(\tau)\Big\|_{\dot{H}^{[\frac{d}{p}]-1}_{\mathcal{L}^{\tilde p, \infty}}}
\Big\|v(\tau)\Big\|_{\dot{H}^{[\frac{d}{p}]-1}_{\mathcal{L}^{\tilde p, \infty}}} \\
=(t - \tau)^{\frac{\alpha}{2}-1}\Big\|u(\tau)\Big\|_{\dot{H}^{[\frac{d}{p}]
-1}_{\mathcal{L}^{\tilde p, \infty}}}\Big\|v(\tau)\Big\|_{\dot{H}^{[\frac{d}{p}]
-1}_{\mathcal{L}^{\tilde p, \infty}}},
\end{gather*}
where
$$
\alpha = \alpha\Big(\frac{d}{[\frac{d}{p}]},\tilde p\Big) = \Big[\frac{d}{p}\Big]-\frac{d}{\tilde p}\ ,
$$
this gives the desired result
\begin{gather*}
\Big\|B(u, v)(t)\Big\|_{\dot{H}^{[\frac{d}{p}]-1}_{\mathcal{L}^{\tilde p, 1}}} 
\lesssim  \int_0^t (t - \tau)^{\frac{\alpha}{2}-1}\Big\|u(\tau)
\Big\|_{\dot{H}^{[\frac{d}{p}]-1}_{\mathcal{L}^{\tilde p, \infty}}}
\Big\|v(\tau)\Big\|_{\dot{H}^{[\frac{d}{p}]-1}_{\mathcal{L}^{\tilde p, \infty}}}\dif\tau\notag\\
\leq\int_0^t (t - \tau)^{\frac{\alpha}{2} - 1}\tau^{-\alpha}\underset{0 < \eta < t}
{\rm sup}\eta^{\frac{\alpha}{2}}\Big\|u(\eta)
\Big\|_{\dot{H}^{[\frac{d}{p}]-1}_{\mathcal{L}^{\tilde p, \infty}}}
\underset{0 < \eta < t}{\rm sup}
\eta^{\frac{\alpha}{2}}\Big\|v(\eta)\Big\|_{\dot{H}^{[\frac{d}{p}]
-1}_{\mathcal{L}^{\tilde p, \infty}}}\dif\tau
\end{gather*}
\begin{gather}
= \underset{0 < \eta < t} {\rm sup}\eta^{\frac{\alpha}{2}}\Big\|u(\eta)
\Big\|_{\dot{H}^{[\frac{d}{p}]-1}_{\mathcal{L}^{\tilde p, \infty}}}
\underset{0 < \eta < t}{\rm sup}
\eta^{\frac{\alpha}{2}}\Big\|v(\eta)\Big\|_{\dot{H}^{[\frac{d}{p}]
-1}_{\mathcal{L}^{\tilde p, \infty}}}\int_0^t (t - \tau)^{\frac{\alpha}{2} - 1}
\tau^{-\alpha}\dif\tau \notag \\
\simeq t^{-\frac{\alpha}{2}}\underset{0 < \eta < t}{\rm sup}\eta^{\frac{\alpha}{2}}
\Big\|u(\eta)\Big\|_{\dot{H}^{[\frac{d}{p}]-1}_{\mathcal{L}^{\tilde p, \infty}}}
\underset{0 < \eta < t}{\rm sup}\eta^{\frac{\alpha}{2}}
\Big\|v(\eta)\Big\|_{\dot{H}^{[\frac{d}{p}]-1}_{\mathcal{L}^{\tilde p, \infty}}}. \label{lem.3.5.10}
\end{gather}
Now we check the validity of condition \eqref{3.1.1} for the bilinear term $B(u,v)(t)$. 
From \eqref{lem.3.5.10} we infer that 
$$
\underset{t \rightarrow 0}{\rm lim\ }t^{\frac{\alpha}{2}}
\Big\|B(u, v)(t)\Big\|_{\dot{H}^{[\frac{d}{p}]-1}_{\mathcal{L}^{\tilde p, 1}}}  = 0,
$$
whenever
$$
\underset{t \rightarrow 0}{\rm lim\ }t^{\frac{\alpha}{2}}
\Big\|u(t)\Big\|_{\dot{H}^{[\frac{d}{p}]-1}_{\mathcal{L}^{\tilde p, \infty}}} 
= \underset{t \rightarrow 0}{\rm lim\ }t^{\frac{\alpha}{2}}
\big\|v(t)\Big\|_{\dot{H}^{[\frac{d}{p}]-1}_{\mathcal{L}^{\tilde p, \infty}}} = 0.
$$
Finally, the estimate  \eqref{lem.3.5.2} can be deduced from 
the inequality \eqref{lem.3.5.10}. \qed 
\begin{Th}\label{th.3.1}
Let $1 <p\leq d$ and $1 \leq r < \infty$. Then for all $\tilde p$ be such that
$$
\frac{1}{2p}
+\frac{[\frac{d}{p}]-1}{2d}< \frac{1}{\tilde p} 
< {\rm min}\Big\{\frac{[\frac{d}{p}]}{d},\frac{1}{2}+\frac{[\frac{d}{p}]-1}{2d}\Big\},
$$
there exists a positive constant $\delta_{p, \tilde p,d}$ 
such that for all $T > 0$ and for all $u_0 \in \dot{H}^{\frac{d}{p}
-1}_{\mathcal{L}^{p,r}} (\mathbb{R}^d)$ with ${\rm div}(u_0) = 0$ satisfying
\begin{equation}\label{th.3.1.1}
\underset{0 < t < T}{\rm sup}t^{\frac{1}{2}([\frac{d}{p}] 
- \frac{d}{\tilde p})}\Big\|e^{t\Delta}u_0\Big\|_{\dot{H}^{[\frac{d}{p}]
-1}_{\mathcal{L}^{\tilde p, \infty}}} \leq \delta_{p, \tilde p,d},
\end{equation}
NSE has a unique mild solution $u \in \mathcal{K}^{\tilde p}_{\frac{d}{[\frac{d}{p}]},1,T} \cap L^\infty\big([0, T]; \dot{H}^{\frac{d}{p} - 1}_{\mathcal{L}^{p,r}}\big)$.\\
In particular, the inequality \eqref{th.3.1.1} holds for arbitrary 
$u_0 \in \dot{H}^{\frac{d}{p}-1}_{\mathcal{L}^{p,r}} (\mathbb{R}^d)$ 
when $T(u_0)$ is small enough, and there exists a positive 
constant $\sigma_{p, \tilde p,d}$ such that we can take  $T = \infty$ 
whenever $\Big\|u_0\Big\|_{\dot{B}^{\frac{d}{\tilde p} 
- 1, \infty}_{\mathcal{L}^{\tilde p, \infty}}} \leq \sigma_{p, \tilde p,d}$.
\end{Th}
\textbf{Proof}. From Lemmas \ref{lem3.1} and  \ref{lem.3.5}, the bilinear 
operator $B(u,v)(t)$ is continuous from $\mathcal{K}^{\tilde p}_{\frac{d}{[\frac{d}{p}]},\infty,T} \times\mathcal{K}^{\tilde p}_{\frac{d}{[\frac{d}{p}]},\infty,T}$ into 
$\mathcal{K}^{\tilde p}_{\frac{d}{[\frac{d}{p}]},1,T}$ and we have the inequality
$$
\Big\|B(u, v)\Big\|_{\mathcal{K}^{\tilde p}_{\frac{d}{[\frac{d}{p}]},\infty,T}} 
\leq \Big\|B(u, v)\Big\|_{\mathcal{K}^{\tilde p}_{\frac{d}{[\frac{d}{p}]},1,T}} \leq
 C_{p,\tilde p,d}\Big\|u\Big\|_{\mathcal{K}^{\tilde p}_{\frac{d}{[\frac{d}{p}]},\infty,T}}
\Big\|v\Big\|_{\mathcal{K}^{\tilde p}_{\frac{d}{[\frac{d}{p}]},\infty,T}},
$$
where $C_{p,\tilde p,d}$ is positive constant independent of  $T$. 
From Theorem \ref{th.2.4} and the above inequality, we deduce that 
for any $u_0 \in \dot{H}^{\frac{d}{p}-1}_{\mathcal{L}^{p,r}}$ such that
$$
\Big\|e^{t\Delta}u_0\Big\|_{\mathcal{K}^{\tilde p}_{\frac{d}{[\frac{d}{p}]},\infty,T}} 
= \underset{0 < t < T}{\rm sup}t^{\frac{1}{2}([\frac{d}{p}]
 - \frac{d}{\tilde p})}\Big\|e^{t\Delta}u_0\Big\|_{\dot{H}^{[\frac{d}{p}]
-1}_{\mathcal{L}^{\tilde p, \infty}}} \leq  \frac{1}{4C_{p,\tilde p,d}},
$$
the Navier-Stokes equations has a solution $u$ on the interval $(0, T)$ so that 
\begin{equation}\label{th.3.1.2}
u \in \mathcal{K}^{\tilde p}_{\frac{d}{[\frac{d}{p}]},\infty,T}.
\end{equation}
From Lemmas \ref{lem3.1} and \ref{lem.3.4}, and \eqref{th.3.1.2}, we have 
$$
B(u,u) \in \mathcal{K}^{p}_{p,1,T} \subseteq \mathcal{K}^{p}_{p,r,T}
 \subseteq L^\infty\big([0, T]; \dot{H}^{\frac{d}{p} - 1}_{\mathcal{L}^{p,r}}\big).
$$
From Lemma \ref{lem.2.5}, we also have  
$e^{t\Delta}u_0 \in L^\infty\big([0, T]; \dot{H}^{\frac{d}{p} 
- 1}_{\mathcal{L}^{p,r}}\big)$. Therefore
$$
u = e^{t\Delta}u_0 - B(u,u) \in  L^\infty\big([0, T]; \dot{H}^{\frac{d}{p}
 - 1}_{\mathcal{L}^{p,r}}\big).
$$
For all $u_0 \in \dot{H}^{\frac{d}{p}-1}_{\mathcal{L}^{p,r}}$, 
applying Theorem \ref{th.2.3}, we deduce that 
\begin{equation}\label{th.3.1.3}
u_0 \in \dot{H}^{[\frac{d}{p}]-1}_{\mathcal{L}^{d/[\frac{d}{p}],r}}.
\end{equation}
From \eqref{th.3.1.3}, applying Lemma \ref{lem.3.2}, we get 
$e^{t\Delta}u_0 \in \mathcal{K}^{\tilde p}_{\frac{d}{[\frac{d}{p}]},\infty,T}$. 
From the definition of $\mathcal{K}^{\tilde p}_{p,r,T}$, we deduce that the left-hand 
side of the inequality \eqref{th.3.1.1} converges to $0$ when $T$ 
tends to $0$. Therefore the inequality \eqref{th.3.1.1} holds for arbitrary 
$u_0 \in \dot{H}^{\frac{d}{p}-1}_{\mathcal{L}^{p,r}}$ 
when $T(u_0)$ is small enough. Applying Lemmas \ref{lem.3.2} and \ref{lem.3.5}, 
we conclude that $u \in \mathcal{K}^{\tilde p}_{\frac{d}{[\frac{d}{p}]},1,T}$. \\
Next, applying Theorem 5.4  (\cite{P. G. Lemarie-Rieusset 2002}, p. 45), 
we deduce that the two quantities $\Big\|u_0\Big\|_{\dot{B}^{\frac{d}{\tilde p} 
- 1, \infty}_{\mathcal{L}^{\tilde p, \infty}}}$ and $\underset{0 < t < \infty}{\rm sup}
t^{\frac{1}{2}([\frac{d}{p}] - \frac{d}{\tilde p})}\Big\|e^{t\Delta}
u_0\Big\|_{\dot{H}^{[\frac{d}{p}]-1}_{\mathcal{L}^{\tilde p, \infty}}}$ 
are equivalent, then there \linebreak exists a positive constant 
$\sigma_{p,\tilde p,d}$ such that  $T = \infty$ and  \eqref{th.3.1.1} 
holds whenever $\Big\|u_0\Big\|_{\dot{B}^{\frac{d}{\tilde p} 
- 1, \infty}_{\mathcal{L}^{\tilde p, \infty}}} \leq \sigma_{p, \tilde p,d}$ .\qed 
\begin{Not}
From Theorem \ref{th.2.3} and the proof of Lemma \ref{lem.3.2}, 
and Theorem 5.4  (\cite{P. G. Lemarie-Rieusset 2002}, p. 45), 
we have the following imbedding maps 
$$
\dot{H}^{\frac{d}{p}-1}_{\mathcal{L}^{p,r}}(\mathbb{R}^d)
 \hookrightarrow \dot{H}^{[\frac{d}{p}]-1}_{\mathcal{L}^{d/[\frac{d}{p}],r}}(\mathbb{R}^d) \hookrightarrow \dot{B}^{\frac{d}{\tilde p} - 1, \infty}_{\mathcal{L}^{\tilde p, 1}}(\mathbb{R}^d) \hookrightarrow \dot{B}^{\frac{d}{\tilde p} - 1, \infty}_{\mathcal{L}^{\tilde p, \infty}}(\mathbb{R}^d).
$$
On the other hand, a function in 
$\dot{H}^{\frac{d}{p}-1}_{\mathcal{L}^{p,r}}(\mathbb{R}^d)$ 
can be arbitrarily large in the $\dot{H}^{\frac{d}{p}-1}_{\mathcal{L}^{p,r}}(\mathbb{R}^d)$ 
norm but small in the $\dot{B}^{\frac{d}{\tilde p} - 1, \infty}_{\mathcal{L}^{\tilde p, 
\infty}}(\mathbb{R}^d)$ norm. 
\end{Not}
\vskip 0.4cm
\textbf{3.2. Solutions to the Navier-Stokes equations with the initial value 
in the critical spaces  $\dot{H}^{\frac{d}{p} - 1}_{\mathcal{L}^{p,r}}(\mathbb{R}^d)$ 
with $d \leq p < \infty$ and $1 \leq r < \infty$}.
\begin{Lem}\label{lem.3.2.1} Suppose that 
$u_0 \in \dot{H}^{\frac{d}{p} - 1}_{\mathcal{L}^{p, r}}$ 
with $d \leq p < \infty$ and $1 \leq r < \infty$. 
Then $e^{t\Delta}u_0 \in \mathcal{K}^{\tilde p}_{d,1,\infty}$ for all $\tilde p > p$.
\end{Lem}
\textbf{Proof}.  We prove that
$$
\underset{0 < t < \infty}{\rm sup}t^{\frac{\alpha}{2}}
\big\|e^{t\Delta}u_0\big\|_{\mathcal{L}^{\tilde p, 1}} 
\lesssim \Big\|u_0\Big\|_{\dot{H}^{\frac{d}{p} - 1}_{\mathcal{L}^{p, r}}},
$$
where
$$
\alpha = \alpha (d,\tilde p) = 1-\frac{d}{\tilde p}.
$$
Let $p'$ and $\tilde p'$ be such that
$$
\frac{1}{p} + \frac{1}{p'} = 1, \frac{1}{\tilde p} + \frac{1}{\tilde p'} = 1.
$$
We have
\begin{gather*}
\big\|e^{t\Delta}u_0\big\|_{\mathcal{L}^{\tilde p, 1}} 
= \big\|e^{-t\left| \xi \right|^2}\hat u_0(\xi)\big\|_{L^{\tilde p', 1}_\xi} 
= \big\|e^{-t\left| \xi \right|^2}|\xi|^{1-\frac{d}{p}}|\xi|^{\frac{d}{p} 
- 1}\hat u_0(\xi)\big\|_{L^{\tilde p', 1}_\xi}.
\end{gather*}
Applying Holder's inequality in the Lorentz spaces to obtain
\begin{gather*}
\Big\|e^{-t\left| \xi \right|^2}|\xi|^{1-\frac{d}{p}}|\xi|^{\frac{d}{p} - 1}
\hat u_0(\xi)\Big\|_{L^{\tilde p', 1}_\xi} = \Big\|e^{-t|\xi|^2}|\xi|^{1-\frac{d}{p}}
\Big\|_{L_{\xi}^{\frac{p\tilde p}{\tilde p-p }, 1}}\big\||\xi|^{\frac{d}{p} - 1}
\hat u_0(\xi)\big\|_{L^{p', \infty}}\\
=t^{-\frac{1}{2}(1-\frac{d}{\tilde p})}\Big\|e^{-|\xi|^2}|\xi|^{1-\frac{d}{p}}
\Big\|_{L^{\frac{p\tilde p}{\tilde p-p }, 1}}\big\||\xi|^{\frac{d}{p} - 1}
\hat u_0(\xi)\big\|_{L^{p', \infty}} \\
\simeq t^{-\frac{\alpha}{2}}
\big\||\xi|^{\frac{d}{p} - 1}\hat u_0(\xi)\big\|_{L^{p', \infty}}  
\lesssim t^{-\frac{\alpha}{2}}\big\||\xi|^{\frac{d}{p} - 1}
\hat u_0(\xi)\big\|_{L^{p', r}} = t^{-\frac{\alpha}{2}}
\Big\|u_0\Big\|_{\dot{H}^{\frac{d}{p} - 1}_{\mathcal{L}^{p, r}}}.
\end{gather*}
Therefore this gives the desired result
$$
\big\|e^{t\Delta}u_0\big\|_{\mathcal{L}^{\tilde p, 1}} 
\lesssim t^{-\frac{\alpha}{2}}\Big\|u_0\Big\|_{\dot{H}^{\frac{d}{p} 
- 1}_{\mathcal{L}^{p, r}}}. 
$$
We claim now that
\begin{equation*}
\underset{t \rightarrow 0}{\rm lim\ }t^{\frac{\alpha}{2}}
\big\|e^{t\Delta}u_0\big\|_{\mathcal{L}^{\tilde p, 1}} = 0.
\end{equation*}
For any $\epsilon > 0$. Applying Lemma \ref{lem.3.3} and from the above proof we deduce that
\begin{gather*}
t^{\frac{\alpha}{2}}\big\|e^{t\Delta}u_0\big\|_{\mathcal{L}^{\tilde p, 1}} \leq  \\
t^{\frac{\alpha}{2}}\Big\|e^{-t\left| \xi \right|^2}|\xi|^{1
-\frac{d}{p}}1_{B_n^c}|\xi|^{\frac{d}{p}-1}\hat u_0(\xi)
\Big\|_{L^{\tilde p', 1}_\xi} + 
t^{\frac{\alpha}{2}}\Big\|e^{-t\left| \xi \right|^2}|\xi|^{1-\frac{d}{p}}
1_{B_n}|\xi|^{\frac{d}{p}-1}\hat u_0(\xi)\Big\|_{L^{\tilde p', 1}_\xi} \leq \\
\end{gather*}
\begin{gather*}
C_1\Big\|e^{-|\xi|^2}|\xi|^{1-\frac{d}{p}}
\Big\|_{L^{\frac{p\tilde p}{\tilde p-p }, 1}}
\Big\|1_{B^c_n}|\xi|^{\frac{d}{p} - 1}\hat u_0(\xi)\Big\|_{L^{p', \infty}}\\ 
+\ C_2t^{\frac{\alpha}{2}}\Big\|1_{B_n}|\xi|^{1-\frac{d}{p}}
\Big\|_{L^{\frac{p\tilde p}{\tilde p - p}, 1}}\big\||\xi|^{\frac{d}{p}
-1}\hat u_0(\xi)\big\|_{L^{p', \infty}}\\
 \leq C_3\Big\|1_{B^c_n}|\xi|^{\frac{d}{p} - 1}\hat u_0(\xi)
\Big\|_{L^{p', r}}+ C_4(n)t^{\frac{\alpha}{2}}
\Big\|u_0\Big\|_{\dot{H}^{\frac{d}{p} - 1}_{\mathcal{L}^{p, r}}} < \epsilon
\end{gather*}
for large enough $n$ and small enough $t=t(n) > 0$. \qed
\begin{Lem}\label{th.3.2.2}
Let 
\begin{equation}\label{th.3.2.3}
p \geq d \ \text{and}\ d < \tilde p < 2p.
\end{equation}
Then the bilinear operator $B(u, v)(t)$ is continuous from 
$\mathcal{K}^{\tilde p}_{d,\infty,T} \times \mathcal{K}^{\tilde p}_{d,\infty,T}$ 
into $\mathcal{K}^{p}_{p,1,T}$, and we have the inequality
\begin{equation}\label{th.3.2.4}
\big\|B(u, v)\big\|_{\mathcal{K}^{p}_{p,1,T}} \leq
 C\big\|u\big\|_{\mathcal{K}^{\tilde p}_{d,\infty,T}}
\big\|v\big\|_{\mathcal{K}^{\tilde p}_{d,\infty,T}},
\end{equation}
where C is a positive constant and independent of T.
\end{Lem}
\textbf{Proof}. First, arguing as in Lemma \ref{lem.3.4}, we derive
\begin{gather*}
\dot{\Lambda}^{\frac{d}{p}-1}e^{(t - \tau)\Delta}\mathbb{P}
\nabla .\big(u(\tau) \otimes v(\tau)\big) 
=\frac{1}{(t - \tau)^{\frac{d}{2}(\frac{1}{p}+1)}}
K\Big(\frac{.}{\sqrt{t - \tau}}\Big)*\big(u(\tau) \otimes v(\tau)\big),
\end{gather*}
where the tensor $K(x) = \{K_{l, k, j}(x)\}$ is given by the formula
\begin{gather}
\widehat{K_{l, k, j}}(\xi) = \frac{1}{(2\pi)^{d/2}}|\xi|^{\frac{d}{p}-1}e^{-|\xi|^2}
\Big(\delta_{jk} - \frac{\xi_j\xi_k}{|\xi|^2}\Big)(i\xi_l).\label{th.3.2.6}
\end{gather}
Applying Theorem \ref{th.2.2} for the convolution in the Fourier-Lorentz spaces, we have
\begin{gather}
\Big\|\dot{\Lambda}^{\frac{d}{p}-1}e^{(t - \tau)\Delta}\mathbb{P}\nabla .
\big(u(\tau) \otimes v(\tau)\big)\Big\|_{\mathcal{L}^{p, 1}} \notag \\
\lesssim \frac{1}{(t - \tau)^{\frac{d}{2}(\frac{1}{p}+1)}}
\Big\|K\Big(\frac{.}{\sqrt{t - \tau}}\Big)\Big\|_{\mathcal{L}^{r, 1}}
\big\|\big(u(\tau) \otimes v(\tau)\big)\big\|_{\mathcal{L}^{\frac{\tilde p}{2}, \infty}},\label{th.3.2.7}
\end{gather}
where
\begin{gather}
\frac{1}{r} = 1+ \frac{1}{p} - \frac{2}{\tilde p}.\label{th.3.2.8}
\end{gather}
Note that from the inequality \eqref{th.3.2.3}, we can check that $1 < r < \infty$. 
Applying Theorem \ref{th.2.1}, we have
\begin{gather}
\big\|u(\tau) \otimes v(\tau)\big\|_{\mathcal{L}^{\frac{\tilde p}{2}, \infty}} 
\lesssim \big\|u(\tau)\big\|_{\mathcal{L}^{\tilde p, \infty}}
\big\|v(\tau)\big\|_{\mathcal{L}^{\tilde p, \infty}}.
\label{th.3.2.9}
\end{gather}
From the equalities \eqref{th.3.2.6} and \eqref{th.3.2.8} it follows that
\begin{gather}
\Big\|K\Big(\frac{.}{\sqrt{t - \tau}}\Big)\Big\|_{\mathcal{L}^{r,1}} 
=  (t - \tau)^{\frac{d}{2r}}\big\|\hat K\big\|_{L^{r', 1}} 
\simeq (t - \tau)^{\frac{d}{2}(1+ \frac{1}{p} - \frac{2}{\tilde p})}.
\label{th.3.2.10}
\end{gather}
From the estimates \eqref{th.3.2.7}, \eqref{th.3.2.9}, and \eqref{th.3.2.10}, we deduce that 
\begin{gather*}
\Big\|e^{(t - \tau)\Delta}\mathbb{P}\nabla .\big(u(\tau) 
\otimes v(\tau)\big)\Big\|_{\dot{H}^{\frac{d}{p} - 1}_{\mathcal{L}^{p, 1}}}
\lesssim (t - \tau)^{-\frac{d}{\tilde p}}\big\|u(\tau)
\big\|_{\mathcal{L}^{\tilde p, \infty}}\big\|v(\tau)\big\|_{\mathcal{L}^{\tilde p, \infty}}\\
 = (t - \tau)^{\alpha-1}\big\|u(\tau)\big\|_{\mathcal{L}^{\tilde p, \infty}}
\big\|v(\tau)\big\|_{\mathcal{L}^{\tilde p, \infty}},
\end{gather*}
where
$$
\alpha = \alpha (d,\tilde p) = 1-\frac{d}{\tilde p}.
$$
This gives the desired result
\begin{gather}
\Big\|B(u, v)(t)\Big\|_{\dot{H}^{\frac{d}{p} - 1}_{\mathcal{L}^{p, 1}}} 
\lesssim  \int_0^t (t - \tau)^{\alpha -1}\big\|u(\tau)\big\|_{\mathcal{L}^{\tilde p, \infty}}
\big\|v(\tau)\big\|_{\mathcal{L}^{\tilde p, \infty}}\dif\tau \notag \\
\leq\int_0^t (t - \tau)^{\alpha-1}\tau^{-\alpha}\underset{0 < \eta < t}{\rm sup}
\eta^{\frac{\alpha}{2}}\big\|u(\eta)\big\|_{\mathcal{L}^{\tilde p, \infty}}
\underset{0 < \eta < t}{\rm sup}\eta^{\frac{\alpha}{2}}
\big\|v(\eta)\big\|_{\mathcal{L}^{\tilde p, \infty}}\dif\tau \notag \\
= \underset{0 < \eta < t}{\rm sup}\eta^{\frac{\alpha}{2}}
\big\|u(\eta)\big\|_{\mathcal{L}^{\tilde p, \infty}}
\underset{0 < \eta < t}{\rm sup}\eta^{\frac{\alpha}{2}}
\big\|v(\eta)\big\|_{\mathcal{L}^{\tilde p, \infty}}\int_0^t (t - \tau)^{\alpha-1
}\tau^{-\alpha}\dif\tau \notag \\
\simeq \underset{0 < \eta < t}{\rm sup}\eta^{\frac{\alpha}{2}}
\big\|u(\eta)\big\|_{\mathcal{L}^{\tilde p, \infty}}
\underset{0 < \eta < t}{\rm sup}\eta^{\frac{\alpha}{2}}
\big\|v(\eta)\big\|_{\mathcal{L}^{\tilde p, \infty}} . \label{th.3.2.11}
\end{gather}
From \eqref{th.3.2.11} it follows the validity of \eqref{3.1.2}  since  
$$
\underset{t \rightarrow 0}{\rm lim}\Big\|B(u,v)(t)\Big\|_{\dot{H}^{\frac{d}{p} 
- 1}_{\mathcal{L}^{p, 1}}} = 0,
$$
whenever
$$
\underset{t \rightarrow 0}{\rm lim\ }t^{\frac{\alpha}{2}}
\big\|u(t)\big\|_{\mathcal{L}^{\tilde p, \infty}} 
= \underset{t \rightarrow 0}{\rm lim\ }t^{\frac{\alpha}{2}}
\big\|v(t)\big\|_{\mathcal{L}^{\tilde p, \infty}} = 0.
$$
The estimate  \eqref{th.3.2.4} can be deduced from the inequality \eqref{th.3.2.11}. \qed 

\begin{Lem}\label{lem.3.2.3}
Let $\tilde p > d$, then the bilinear operator $B(u, v)(t)$ is continuous from 
$\mathcal{K}^{\tilde p}_{d,\infty,T} \times \mathcal{K}^{\tilde p}_{d,\infty,T}$ 
into $\mathcal{K}^{\tilde p}_{d,1,T}$, and we have the inequality
\begin{equation}\label{lem.3.2.3.1}
\big\|B(u, v)\big\|_{\mathcal{K}^{\tilde p}_{d,1,T}} \leq
 C\big\|u\big\|_{\mathcal{K}^{\tilde p}_{d,\infty,T}}
\big\|v\big\|_{\mathcal{K}^{\tilde p}_{d,\infty,T}},
\end{equation}
where C is a positive constant and independent of T.
\end{Lem}
\textbf{Proof}. First, arguing as in Lemma \ref{lem.3.4}, we derive
\begin{gather*}
e^{(t - \tau)\Delta}\mathbb{P}\nabla .\big(u(\tau) \otimes v(\tau)\big) 
=\frac{1}{(t - \tau)^{\frac{d+1}{2}}}K\Big(\frac{.}{\sqrt{t - \tau}}\Big)*\big(u(\tau)
 \otimes v(\tau)\big),
\end{gather*}
where the tensor $K(x) = \{K_{l, k, j}(x)\}$ is given by the formula
\begin{gather}
\widehat{K_{l, k, j}}(\xi) = \frac{1}{(2\pi)^{d/2}}e^{-|\xi|^2}
\Big(\delta_{jk} - \frac{\xi_j\xi_k}{|\xi|^2}\Big)(i\xi_l).\label{lem.3.2.3.3}
\end{gather}
Applying Theorem \ref{th.2.2} for the convolution in the Fourier-Lorentz spaces, we have
\begin{gather}
\Big\|e^{(t - \tau)\Delta}\mathbb{P}\nabla .
\big(u(\tau) \otimes v(\tau)\big)\Big\|_{\mathcal{L}^{\tilde p, 1}} \notag \\
\lesssim \frac{1}{(t - \tau)^{\frac{d+1}{2}}}
\Big\|K\Big(\frac{.}{\sqrt{t - \tau}}\Big)\Big\|_{\mathcal{L}^{r, 1}}
\big\|\big(u(\tau) \otimes v(\tau)\big)
\big\|_{\mathcal{L}^{\frac{\tilde p}{2}, \infty}},\label{lem.3.2.3.4}
\end{gather}
where
\begin{gather}
\frac{1}{r} = 1 - \frac{1}{\tilde p}.\label{lem.3.2.3.5}
\end{gather}
Applying Theorem \ref{th.2.1}, we have
\begin{gather}
\big\|u(\tau) \otimes v(\tau)\big\|_{\mathcal{L}^{\frac{\tilde p}{2}, 
\infty}} \lesssim \big\|u(\tau)\big\|_{\mathcal{L}^{\tilde p, \infty}}
\big\|v(\tau)\big\|_{\mathcal{L}^{\tilde p, \infty}}.
\label{lem.3.2.3.6}
\end{gather}
From the equalities \eqref{lem.3.2.3.3} and \eqref{lem.3.2.3.5} it follows that
\begin{gather}
\Big\|K\Big(\frac{.}{\sqrt{t - \tau}}\Big)\Big\|_{\mathcal{L}^{r,1}} 
= (t - \tau)^{\frac{d}{2r}}\big\|\hat K\big\|_{L^{r', 1}} 
\simeq (t - \tau)^{\frac{d}{2}(1 - \frac{1}{\tilde p})}.
\label{lem.3.2.3.7}
\end{gather}
From the estimates \eqref{lem.3.2.3.4}, \eqref{lem.3.2.3.6}, 
and \eqref{lem.3.2.3.7}, we deduce that 
\begin{gather*}
\Big\|e^{(t - \tau)\Delta}\mathbb{P}\nabla .
\big(u(\tau) \otimes v(\tau)\big)\Big\|_{\mathcal{L}^{\tilde p, 1}} 
\lesssim (t - \tau)^{ - \frac{1}{2}(\frac{d}{\tilde p}+1)}
\big\|u(\tau)\big\|_{\mathcal{L}^{\tilde p, \infty}}
\big\|v(\tau)\big\|_{\mathcal{L}^{\tilde p, \infty}}\\
= (t - \tau)^{\frac{\alpha}{2}-1}\big\|u(\tau)\big\|_{\mathcal{L}^{\tilde p, \infty}}
\big\|v(\tau)\big\|_{\mathcal{L}^{\tilde p, \infty}},
\end{gather*}
where
$$
\alpha = \alpha (d,\tilde p) = 1-\frac{d}{\tilde p}.
$$
This gives the desired result
\begin{gather}
\big\|B(u, v)(t)\big\|_{\mathcal{L}^{\tilde p, 1}} \lesssim  
\int_0^t (t - \tau)^{\frac{\alpha}{2}-1}\big\|u(\tau)\big\|_{\mathcal{L}^{\tilde p, 
\infty}}\big\|v(\tau)\big\|_{\mathcal{L}^{\tilde p, \infty}}\dif\tau\notag \\
 \leq\int_0^t (t - \tau)^{\frac{\alpha}{2}-1}\tau^{-\alpha}
\underset{0 < \eta < t}{\rm sup}\eta^{\frac{\alpha}{2}}
\big\|u(\eta)\big\|_{\mathcal{L}^{\tilde p, \infty}}\underset{0 < \eta < t}{\rm sup}
\eta^{\frac{\alpha}{2}}\big\|v(\eta)\big\|_{\mathcal{L}^{\tilde p, \infty}}\dif\tau \notag \\
= \underset{0 < \eta < t}{\rm sup}\eta^{\frac{\alpha}{2}}
\big\|u(\eta)\big\|_{\mathcal{L}^{\tilde p, \infty}}
\underset{0 < \eta < t}{\rm sup}\eta^{\frac{\alpha}{2}}
\big\|v(\eta)\big\|_{\mathcal{L}^{\tilde p, \infty}}\int_0^t (t - \tau)^{\frac{\alpha}{2}
-1}\tau^{-\alpha}\dif\tau \notag \\
\simeq t^{-\frac{\alpha}{2}}\underset{0 < \eta < t}{\rm sup}
\eta^{\frac{\alpha}{2}}\big\|u(\eta)\big\|_{\mathcal{L}^{\tilde p, \infty}}
\underset{0 < \eta < t}{\rm sup}\eta^{\frac{\alpha}{2}}
\big\|v(\eta)\big\|_{\mathcal{L}^{\tilde p, \infty}}. \label{lem.3.2.3.8}
\end{gather}
From \eqref{lem.3.2.3.8} it follows the validity \eqref{3.1.1} since 
$$
\underset{t \rightarrow 0}{\rm lim\ }t^{\frac{\alpha}{2}}
\big\|B(u,v)(t)\big\|_{\mathcal{L}^{\tilde p, 1}} = 0,
$$
whenever
$$
\underset{t \rightarrow 0}{\rm lim\ }t^{\frac{\alpha}{2}}
\big\|u(t)\big\|_{\mathcal{L}^{\tilde p, \infty}} 
= \underset{t \rightarrow 0}{\rm lim\ }t^{\frac{\alpha}{2}}
\big\|v(t)\big\|_{\mathcal{L}^{\tilde p, \infty}} = 0.
$$
Finally, the estimate  \eqref{lem.3.2.3.1} can be deduced 
from the inequality \eqref{lem.3.2.3.8}. \qed \\
The following lemma is a generalization of Lemma \ref{lem.3.2.3}.
\begin{Lem}\label{lem.3.2.4} Let $d < \tilde p_1 < \infty$ 
and $d \leq \tilde p_2 < \infty $ be such that one of the following conditions is satisfied 
$$
d < \tilde p_1 < 2d, d \leq \tilde p_2 < \frac{d\tilde p_1}{2d - \tilde p_1}, 
$$
or
$$
\tilde p_1 = 2d, d \leq \tilde p_2 < \infty,
$$
or
$$
2d < \tilde p_1< \infty, \frac{\tilde p_1}{2} < \tilde p_2 < \infty. 
$$
Then the bilinear operator $B(u, v)(t)$ is continuous 
from $\mathcal{K}^{\tilde p_1}_{d,\infty,T} \times \mathcal{K}^{\tilde p_1}_{d,\infty,T}$ 
into $\mathcal{K}^{\tilde p_2}_{d,1,T}$, and we have the inequality
\begin{equation*}
\big\|B(u, v)\big\|_{\mathcal{K}^{\tilde p_2}_{d,1,T}} \leq
 C\big\|u\big\|_{\mathcal{K}^{\tilde p_1}_{d,\infty,T}}
\big\|v\big\|_{\mathcal{K}^{\tilde p_1}_{d,\infty,T}},
\end{equation*}
where C is a positive constant and independent of T.
\end{Lem}
\begin{Th}\label{th.3.2.1}
Let $ p \geq d \ and\ 1 \leq r  < \infty$. 
Then for any $\tilde p$ such that
\begin{equation}\label{Co.3.2.1.1}
 \tilde p > p, 
\end{equation}
there exists a positive constant $\delta_{\tilde p,d}$ such that for 
all $T > 0$ and for all \linebreak
 $u_0 \in \dot{H}^{\frac{d}{p} - 1}_{\mathcal{L}^{p,r}}(\mathbb{R}^d),
\ with\ {\rm div}(u_0) = 0$ satisfying
\begin{equation}\label{Co.3.2.1.2}
\underset{0 < t < T}{\rm sup}t^{\frac{1}{2}(1 
- \frac{d}{\tilde p})}\big\|e^{t\Delta}u_0\big\|_{\mathcal{L}^{\tilde p, \infty}} \leq \delta_{\tilde p,d},
\end{equation}
NSE has a unique mild solution $u \in \underset{q > p}{\cap} \mathcal{K}^q_{d,1,T} 
\cap L^\infty([0, T]; \dot{H}^{\frac{d}{p} - 1}_{\mathcal{L}^{p,r}})$.\\
In particular, the inequality \eqref{Co.3.2.1.2} holds for arbitrary $u_0 \in \dot{H}^{\frac{d}{p} - 1}_{\mathcal{L}^{p,r}}(\mathbb{R}^d) $ with $T(u_0)$ is small enough, 
and there exists a positive constant $\sigma_{\tilde p,d}$ such that 
we can take  $T = \infty$ whenever $\Big\|u_0\Big\|_{\dot{B}^{\frac{d}{\tilde p} - 1, \infty}_{\mathcal{L}^{\tilde p,\infty}}} \leq \sigma_{\tilde p,d}$.
\end{Th}
\textbf{Proof}. Applying Lemma  \ref{lem.3.2.3} and 
Theorem \ref{th.2.4}, we deduce that there exists a positive 
constant $\delta_{\tilde p,d}$ such that for all $T > 0$ 
and for all  $u_0 \in \dot{H}^{\frac{d}{p} - 1}_{\mathcal{L}^{p,r}}
(\mathbb{R}^d)$, ${\rm with\ div}(u_0) = 0$ satisfying the 
inequality \eqref{Co.3.2.1.2} then NSE has a unique mild 
solution $u \in \mathcal{K}^{\tilde p}_{d,1,T}$. 
Next, we prove that  $u \in \underset{q > p}{\cap} \mathcal{K}^q_{d,1,T}$.\\
Consider two cases $d < \tilde p < 2d$ and $2d \leq \tilde p < \infty$ separately.\\
First, we consider the case $d < \tilde p < 2d$. 
We consider two possibilities $\tilde p > \frac{4d}{3}$ 
and $\tilde p \leq \frac{4d}{3}$. In the case $\tilde p > \frac{4d}{3}$, 
we apply  Lemmas \ref{lem.3.2.1} and \ref{lem.3.2.4} to 
obtained $u \in \mathcal{K}^q_{d,1,T}$ for all $q$ 
satisfying $p < q < \tilde p_1 $ where $\tilde p_1 = \frac{d\tilde p}{2d - \tilde p} > 2d$. 
Thus, $u \in \mathcal{K}^{2d}_{d,1,T}$. Applying again 
Lemmas \ref{lem.3.2.1} and \ref{lem.3.2.4}, we deduce 
that $u \in \mathcal{K}^q_{d,1,T}$ for all $q > p$. In the case 
$\tilde p \leq \frac{4d}{3}$, we set up the following series of numbers  
$\{\tilde p_i\}_{0 \leq i \leq N}$ by inductive. Set $\tilde p_0 = \tilde p$ 
and $\tilde p_1 = \frac{d\tilde p_0}{2d - \tilde p_0}$. 
We have $\tilde p_1 >  \tilde p_0$. If $\tilde p_1 > \frac{4d}{3}$ 
then set $N=1$ and stop here. In the case $\tilde p_1 \leq \frac{4d}{3}$ 
set $\tilde p_2 = \frac{d\tilde p_1}{2d - \tilde p_1}$. 
We have $\tilde p_2 >  \tilde p_1$. If $\tilde p_2 > \frac{4d}{3}$ then set 
$N=2$ and stop here. In the case $\tilde p_2 \leq \frac{4d}{3}$,   
set $\tilde p_3 = \frac{d\tilde p_2}{2d - \tilde p_2}$.  
We have $\tilde p_3 > \tilde  p_2$, and so on, there exists  $k \geq 0$ 
such that $\tilde p_k \leq \frac{4d}{3}, \tilde p_{k+1}   
= \frac{d\tilde p_k}{2d - \tilde p_k} > \frac{4d}{3}$. 
We set $N = k+1$ and stop here, and we have
\begin{gather*}
\tilde p_0 = \tilde p, \tilde p_i = \frac{d\tilde p_{i-1}}{2d 
- \tilde p_{i-1}}, \tilde p_i > \tilde p_{i-1}\ {\rm for}\ i =1,2,3,..,N,\\
2d \geq \tilde p_{N} > \frac{4d}{3} \geq \tilde p_{N-1}.
\end{gather*}
From  $u \in \mathcal{K}^{\tilde p_0}_{d,1,T}$, applying Lemmas \ref{lem.3.2.1} and \ref{lem.3.2.4} 
to obtained $u \in \mathcal{K}^q_{d,1,T}$ for all $q$ 
satisfying $p < q < \tilde p_1$.  Then applying again 
Lemmas \ref{lem.3.2.1} and \ref{lem.3.2.4} to obtained 
$u \in \mathcal{K}^q_{d,1,T}$ for all $q$ satisfying 
$p < q < \tilde p_2$, and so on, finishing we have 
$u \in \mathcal{K}^q_{d,1,T}$ for all $q$ satisfying 
$p < q < \tilde p_N$.  Therefore $u \in \mathcal{K}^q_{d,1,T}$ for all $q$ satisfying 
$\frac{4d}{3} < q < \tilde p_N$. From the proof  of 
the case $\tilde p > \frac{4d}{3}$, we have $u \in \mathcal{K}^q_{d,1,T}$ for all $q > p$.\\
Next, we consider the case $2d \leq \tilde p < \infty$. Let $i \in \mathbb{N}$ be such that 
$$
\frac{\tilde p}{2^{i-1}} \geq {\rm max}\{2d, p\} > \frac{\tilde p}{2^i}.
$$
From $\tilde p \geq {\rm max}\{2d, p\}$, we have $i \geq 1$. 
Applying the Lemmas \ref{lem.3.2.1} and \ref{lem.3.2.4} 
to obtained $u \in \mathcal{K}^q_{d,1,T}$ for all $ q > \frac{\tilde p}{2}$. 
Applying again Lemmas \ref{lem.3.2.1} and \ref{lem.3.2.4} 
to obtained $u \in \mathcal{K}^q_{d,1,T}$ for all $ q > \frac{\tilde p}{2^2}$, 
and so on, finishing we have $u \in \mathcal{K}^q_{d,1,T}$ 
for all $q > \frac{\tilde p}{2^{i-1}}$. Applying again 
Lemmas \ref{lem.3.2.1} and \ref{lem.3.2.4} to obtained 
$u \in \mathcal{K}^q_{d,1,T}$ for all $q > {\rm max}\{p, \frac{\tilde p}{2^i}\}$. 
If $p \geq \frac{\tilde p}{2^i}$ then we have $u \in \mathcal{K}^q_{d,1,T}$ 
for all $q > p$. If $p < \frac{\tilde p}{2^i}$ then $2d > \frac{\tilde p}{2^i}$. 
Thus $u \in \mathcal{K}^q_{d,1,T}$ for all $q$ satisfying $\frac{\tilde p}{2^i} < q < 2d$. 
Therefore, from the proof of the case $d < \tilde p < 2d$, we 
have $u \in \mathcal{K}^q_{d,1,T}$ for all $q > p$.\\
The fact that $u \in L^\infty([0, T]; \dot{H}^{\frac{d}{p} - 1}_{\mathcal{L}^{p,r}})$ 
can be deduced from Lemmas \ref{lem.2.5}  and \ref{th.3.2.2}. Applying Lemma \ref{lem.3.2.1}, we get 
$e^{t\Delta}u_0 \in \mathcal{K}^{\tilde p}_{d,\infty,T}$. 
From the definition of $\mathcal{K}^{\tilde p}_{p,r,T}$, we deduce that the left-hand 
side of the inequality \eqref{Co.3.2.1.2} converges to $0$ when $T$ 
tends to $0$. Therefore the inequality \eqref{Co.3.2.1.2} holds for arbitrary 
$u_0 \in \dot{H}^{\frac{d}{p}-1}_{\mathcal{L}^{p,r}}$ 
when $T(u_0)$ is small enough. \\
Next, applying Theorem 5.4  (\cite{P. G. Lemarie-Rieusset 2002}, p. 45), 
we deduce that the two quantities $\Big\|u_0\Big\|_{\dot{B}^{\frac{d}{\tilde p} 
- 1, \infty}_{\mathcal{L}^{\tilde p, \infty}}}$ and $\underset{0 < t < \infty}{\rm sup}t^{\frac{1}{2}(1 
- \frac{d}{\tilde p})}\big\|e^{t\Delta}u_0\big\|_{\mathcal{L}^{\tilde p, \infty}}$ 
are equivalent, then there \linebreak exists a positive constant 
$\sigma_{\tilde p,d}$ such that  $T = \infty$ and  \eqref{Co.3.2.1.2} 
holds whenever $\Big\|u_0\Big\|_{\dot{B}^{\frac{d}{\tilde p} 
- 1, \infty}_{\mathcal{L}^{\tilde p, \infty}}} \leq \sigma_{\tilde p,d}$ .\qed 
\begin{Not}
From the proof of 
Lemma \ref{lem.3.2.1} and Theorem 5.4  (\cite{P. G. Lemarie-Rieusset 2002}, p. 45), 
we have the following imbedding maps 
$$
\dot{H}^{\frac{d}{p}-1}_{\mathcal{L}^{p,r}}(\mathbb{R}^d) 
\hookrightarrow \dot{B}^{\frac{d}{\tilde p} - 1, \infty}_{\mathcal{L}^{\tilde p, 
1}}(\mathbb{R}^d) \hookrightarrow \dot{B}^{\frac{d}{\tilde p} - 1, 
\infty}_{\mathcal{L}^{\tilde p, \infty}}(\mathbb{R}^d).
$$
On the other hand, a function in 
$\dot{H}^{\frac{d}{p}-1}_{\mathcal{L}^{p,r}}(\mathbb{R}^d)$ 
can be arbitrarily large in the $\dot{H}^{\frac{d}{p}-1}_{\mathcal{L}^{p,r}}(\mathbb{R}^d)$ 
norm but small in the $\dot{B}^{\frac{d}{\tilde p} - 1, \infty}_{\mathcal{L}^{\tilde p, 
\infty}}(\mathbb{R}^d)$ norm. 
\end{Not}
\vskip 0.8cm
\textbf{3.3. Solutions to the Navier-Stokes equations 
with initial value in the critical spaces  $\dot{H}^{d - 1}_{\mathcal{L}^{1,r}}
(\mathbb{R}^d)$ with $1 \leq r < \infty$}.\\

We define an auxiliary space $\mathcal{K}_{s,r,T}$ which is made 
up by the functions $u(t,x)$ such that 
$$
\big\|u\big\|_{\mathcal{K}_{s,r,T}}:= 
\underset{0 < t < T}{{\rm sup}}t^{\frac{\alpha}{2}}
\Big\|u(t,x)\Big\|_{\dot{H}^s_{\mathcal{L}^{1,r}}} < \infty,
$$
and
\begin{equation}\label{4.1}
\underset{t \rightarrow 0}{\rm lim\ }t^{\frac{\alpha}{2}}
\Big\|u(t,x)\Big\|_{\dot{H}^s_{\mathcal{L}^{1,r}}} = 0,
\end{equation}
with 
$$
d-1 \leq s < d, 1 \leq r \leq \infty, T >0,
$$
and 
$$
\alpha = \alpha(s) =s+1-d.
$$
In the case $s = d-1$, it is also convenient to define the 
space $\mathcal{K}_{d-1,r,T}$ as the natural space  $L^\infty\big([0, T]; \dot{H}^{d-1}_{\mathcal{L}^{1,r}}\big)$ 
with the additional condition that its elements $u(t,x)$ satisfy 
\begin{equation}\label{4.2}
\underset{t \rightarrow 0}{\text{lim}\ }\big\|u(t,x)\big\|_{\dot{H}^{d-1}_{\mathcal{L}^{1,r}}} = 0.
\end{equation}
\begin{Lem}\label{lem4.1} Let $1 \leq r \leq \tilde{r}\leq \infty$. Then we have the following imbedding
$$
\mathcal{K}_{s,1,T} \hookrightarrow \mathcal{K}_{s,r,T} \hookrightarrow \mathcal{K}_{s,\tilde r,T} \hookrightarrow \mathcal{K}_{s,\infty,T}.
$$
\end{Lem}
\textbf{Proof}. It is deduced from Lemma \ref{lem2.1} (a) and 
the definition of  $\mathcal{K}_{s,r,T}$. \qed
\begin{Lem}\label{lem.4.2} Suppose that $u_0 \in \dot{H}^{d-1}_{\mathcal{L}^{1,r}}(\mathbb{R}^d)$ with $1 \leq r < \infty$, then $e^{t\Delta}u_0 \in \mathcal{K}_{s,r,\infty}$ with $d-1 < s <d$.
\end{Lem}
\textbf{Proof}. We prove that
\begin{equation}\label{4.10}
\underset{0 < t < \infty}{\rm sup}t^{\frac{\alpha}{2}}\Big\|e^{t\Delta}
u_0\Big\|_{\dot{H}^s_{\mathcal{L}^{1,r}}} \lesssim 
\Big\|u_0\Big\|_{\dot{H}^{d-1}_{\mathcal{L}^{1,r}}}\ {\rm for}\ 1\leq r \leq \infty.
\end{equation}
We have
\begin{gather}
\Big\|e^{t\Delta}u_0\Big\|_{\dot{H}^s_{\mathcal{L}^{1,r}}} =
\big\|e^{-t\left| \xi \right|^2}|\xi|^s\hat u_0(\xi)\big\|_{L^{\infty,r}_\xi} =
\big\||\xi|^{s+1-d}e^{-t\left| \xi \right|^2}|\xi|^{d-1}
\hat u_0(\xi)\big\|_{L^{\infty, r}_\xi} \notag \\
\leq t^{-\frac{s+1-d}{2}}\big\||\xi|^{s+1-d}e^{-\left|\xi \right|^2}\big\|_{L^{\infty}}
\big\||\xi|^{d-1}\hat u_0(\xi)\big\|_{L^{\infty, r}} \simeq t^{-\frac{\alpha}{2}}
\Big\|u_0\Big\|_{\dot{H}^{d-1}_{\mathcal{L}^{1,r}}}. \label{4.3}
\end{gather}
We claim now that
$$
\underset{t \rightarrow 0}{\rm lim\ }t^{\frac{\alpha}{2}}
\Big\|e^{t\Delta}u_0\Big\|_{\dot{H}^s_{\mathcal{L}^{1,r}}} = 0\ \ {\rm for}\ 1\leq r < \infty. 
$$
From the inequality \eqref{4.3}, we have
\begin{gather*}
t^{\frac{\alpha}{2}}\Big\|e^{t\Delta}u_0\Big\|_{\dot{H}^s_{\mathcal{L}^{1,r}}}
 \leq \\ t^{\frac{\alpha}{2}}\big\||\xi|^{s+1-d}e^{-t\left| \xi \right|^2}1_{B_n^c}|\xi|^{d-1}
\hat u_0(\xi)\big\|_{L^{\infty, r}_\xi} + t^{\frac{\alpha}{2}}\big\||\xi|^{s+1-d}e^{-t\left| \xi \right|^2}1_{B_n}|\xi|^{d-1}
\hat u_0(\xi)\big\|_{L^{\infty, r}_\xi}.
\end{gather*}
For any $\epsilon > 0$, applying Lemma \ref{lem.3.3}, we have
\begin{gather}
t^{\frac{\alpha}{2}}\big\||\xi|^{s+1-d}e^{-t\left| \xi \right|^2}1_{B^c_n}|\xi|^{d-1}
\hat u_0(\xi)\big\|_{L^{\infty, r}_\xi}
\leq \big\||\xi|^{s+1-d}e^{-\left|\xi \right|^2}\big\|_{L^{\infty}}
\big\|1_{B^c_n}|\xi|^{d-1}\hat u_0(\xi)\big\|_{L^{\infty, r}} \notag\\
= C\big\|1_{B^c_n}|\xi|^{d-1}\hat u_0(\xi)\big\|_{L^{\infty, r}} 
< \frac{\epsilon}{2}, \label{lem.4.4}
\end{gather}
for large enough $n$. Fixed one of such $n$, 
we have the following estimates
\begin{gather}
t^{\frac{\alpha}{2}}\big\||\xi|^{s+1-d}e^{-t\left| \xi \right|^2}1_{B_n}
|\xi|^{d-1}\hat u_0(\xi)\big\|_{L^{\infty, r}_\xi} \notag\\
\leq t^{\frac{\alpha}{2}}\big\|1_{B_n}|\xi|^{s+1-d}
e^{-t\left|\xi\right|^2}\big\|_{L^{\infty}}
\big\||\xi|^{d-1}\hat u_0(\xi)\big\|_{L^{\infty, r}} \notag \\
 \leq t^{\frac{\alpha}{2}}\big\|1_{B_n}|\xi|^{s+1-d}\big\|_{L^{\infty}}
\big\||\xi|^{d-1}\hat u_0(\xi)\big\|_{L^{\infty, r}}
= t^{\frac{\alpha}{2}}n^{s+1-d}\big\||\xi|^{d-1}\hat u_0(\xi)\big\|_{L^{\infty, r}} \notag \\
=  t^{\frac{\alpha}{2}}n^{s+1-d}\Big\|u_0\Big\|_{\dot{H}^{d-1}_{\mathcal{L}^{1, r}}} 
< \frac{\epsilon}{2} \label{lem.4.5}
\end{gather}
for small enough $t =t(n) > 0$. From the estimates \eqref{lem.4.4} 
and \eqref{lem.4.5}, we have,
$$
t^{\frac{\alpha}{2}}\Big\|e^{t\Delta}u_0\Big\|_{\dot{H}^s_{\mathcal{L}^{1,r}}} 
\leq  C\big\|1_{B^c_n}|\xi|^{d-1}\hat u_0(\xi)\big\|_{L^{\infty, r}} + 
t^{\frac{\alpha}{2}}n^{s+1-d}\Big\|u_0\Big\|_{\dot{H}^{d-1}_{\mathcal{L}^{1, r}}} 
< \epsilon. \qed
$$
\begin{Lem}\label{lem.4.3}
Let $d-1 <s <d$. Then the bilinear operator $B(u, v)(t)$ is continuous from 
$\mathcal{K}_{s,\infty,T} \times \mathcal{K}_{s,\infty,T}$ 
into $\mathcal{K}_{s,1,T}$ and we have the inequality
\begin{equation}\label{4.5}
\big\|B(u, v)\big\|_{\mathcal{K}_{s,1,T}} \leq
 C\Big\|u\Big\|_{\mathcal{K}_{s,\infty,T}}
\Big\|v\Big\|_{\mathcal{K}_{s,\infty,T}},
\end{equation}
where C is a positive constant and independent of T.
\end{Lem}
\textbf{Proof}.  Using the Fourier transform we get
\begin{gather*}
\mathcal{F}\big(B(u, v)_j(t)\big)(\xi) = \\
\frac{1}{(2\pi)^{\frac{d}{2}}}\int_0^te^{-(t - \tau)|\xi|^2}\sum_{l, k = 1}^d \Big(\delta_{jk} - 
\frac{\xi_j\xi_k}{|\xi|^2}\Big)(i\xi_l)
\big(\widehat{u_l(\tau)}*\widehat{v_k(\tau)}\big)(\xi)\dif\tau.
\end{gather*}
Thus
\begin{gather*}
\big||\xi|^s\mathcal{F}\big(B(u, v)(t)\big)(\xi)\big| \lesssim 
\int_0^t|\xi|^se^{-(t - \tau)|\xi|^2}|\xi|
\big(|\widehat{u(\tau)}|*|\widehat{v(\tau)}|\big)(\xi)\dif\tau.
\end{gather*}
We have
$$
|\xi|^s|\widehat{u(\tau)}(\xi)| \leq \underset{\xi \in \mathbb{R}^d}{\rm sup}
\big||\xi|^s\widehat{u(\tau)}(\xi)\big| = \big\|u(\tau)\big\|_{\dot{H}^s_{\mathcal{L}^1}}
\ {\rm and}\ |\xi|^s|\widehat{v(\tau)}(\xi)| \leq \big\|v(\tau)\big\|_{\dot{H}^s_{\mathcal{L}^1}},
$$
therefore
$$
|\widehat{u(\tau)}(\xi)| \leq \frac{\big\|u(\tau)\big\|_{\dot{H}^s_{\mathcal{L}^1}}}{|\xi|^s},
\ |\widehat{v(\tau)}(\xi)| \leq \frac{\big\|v(\tau)\big\|_{\dot{H}^s_{\mathcal{L}^1}}}{|\xi|^s}.
$$
A standard argument shows that
$$
\frac{1}{|\xi|^s}*\frac{1}{|\xi|^s} = \frac{C}{|\xi|^{2s-d}}.
$$
From the above estimates and Lemma \ref{lem2.1} (b), we have
\begin{gather*}
\big(|\widehat{u(\tau)}|*|\widehat{v(\tau)}|\big)(\xi)
\leq \frac{\big\|u(\tau)\big\|_{\dot{H}^s_{\mathcal{L}^1}}}{|\xi|^s}*
\frac{\big\|v(\tau)\big\|_{\dot{H}^s_{\mathcal{L}^1}}}{|\xi|^s} \simeq\\
 \frac{ \big\|u(\tau)\big\|_{\dot{H}^s_{\mathcal{L}^1}}
\big\|v(\tau)\big\|_{\dot{H}^s_{\mathcal{L}^1}}}{|\xi|^{2s-d}} =
 \frac{ \big\|u(\tau)\big\|_{\dot{H}^s_{\mathcal{L}^{1,\infty}}}
\big\|v(\tau)\big\|_{\dot{H}^s_{\mathcal{L}^{1,\infty}}}}{|\xi|^{2s-d}},
\end{gather*}
this gives the desired result
\begin{gather*}
\int_0^t|\xi|^se^{-(t - \tau)|\xi|^2}|\xi|
\big(|\widehat{u(\tau)}|*|\widehat{v(\tau)}|\big)(\xi)\dif\tau \\ \lesssim
\int_0^t|\xi|^{d+1-s}e^{-(t - \tau)|\xi|^2}
\big\|u(\tau)\big\|_{\dot{H}^s_{\mathcal{L}^{1,\infty}}}\big\|v(\tau)\big\|_{\dot{H}^s_{\mathcal{L}^{1,\infty}}}\dif\tau.
\end{gather*}
Thus
\begin{gather}
\Big\||\xi|^s\mathcal{F}\big(B(u, v)(t)\big)(\xi)\Big\|_{L^{\infty,1}_{\xi}} \lesssim \notag \\
\int_0^t\Big\||\xi|^{d+1-s}e^{-(t - \tau)|\xi|^2}\Big\|_{L^{\infty,1}_{\xi}}
\big\|u(\tau)\big\|_{\dot{H}^s_{\mathcal{L}^{1,\infty}}}
\big\|v(\tau)\big\|_{\dot{H}^s_{\mathcal{L}^{1,\infty}}}\dif\tau  \notag \\
= \int_0^t(t-s)^{\frac{s-d-1}{2}}\Big\||\xi|^{d+1-s}e^{-|\xi|^2}\Big\|_{L^{\infty,1}}
\big\|u(\tau)\big\|_{\dot{H}^s_{\mathcal{L}^{1,\infty}}}
\big\|v(\tau)\big\|_{\dot{H}^s_{\mathcal{L}^{1,\infty}}}\dif\tau   \notag \\
\lesssim \int_0^t (t - \tau)^{\frac{\alpha}{2} - 1}\tau^{-\alpha}\underset{0 < \eta < t}
{\rm sup}\eta^{\frac{\alpha}{2}}\Big\|u(\eta)
\Big\|_{\dot{H}^s_{\mathcal{L}^{1,\infty}}}\underset{0 < \eta < t}{\rm sup}
\eta^{\frac{\alpha}{2}}\Big\|v(\eta)\Big\|_{\dot{H}^s_{\mathcal{L}^{1,\infty}}}\dif\tau \notag \\
= \underset{0 < \eta < t}
{\rm sup}\eta^{\frac{\alpha}{2}}\Big\|u(\eta)
\Big\|_{\dot{H}^s_{\mathcal{L}^{1,\infty}}}\underset{0 < \eta < t}{\rm sup}
\eta^{\frac{\alpha}{2}}\Big\|v(\eta)\Big\|_{\dot{H}^s_{\mathcal{L}^{1,\infty}}}
\int_0^t (t - \tau)^{\frac{\alpha}{2}- 1}\tau^{-\alpha}\dif\tau \notag \\
\simeq t^{-\frac{\alpha}{2}}\underset{0 < \eta < t}{\rm sup}\eta^{\frac{\alpha}{2}}
\Big\|u(\eta)\Big\|_{\dot{H}^s_{\mathcal{L}^{1,\infty}}}
\underset{0 < \eta < t}{\rm sup}\eta^{\frac{\alpha}{2}}
\Big\|v(\eta)\Big\|_{\dot{H}^s_{\mathcal{L}^{1,\infty}}}. \label{4.6}
\end{gather}
Let us now check the validity of the condition \eqref{4.1} 
for the bilinear term $B(u,v)(t)$. 
Indeed, from \eqref{4.6}
$$
\underset{t \rightarrow 0}{\rm lim\ }t^{\frac{\alpha}{2}}
\Big\|B(u, v)(t)\Big\|_{\dot{H}^s_{\mathcal{L}^{1,1}}} 
= \underset{t \rightarrow 0}{\rm lim\ }t^{\frac{\alpha}{2}}
\Big\||\xi|^s\mathcal{F}\big(B(u, v)(t)\big)(\xi)\Big\|_{L^{\infty,1}_{\xi}}  = 0,
$$
whenever
$$
\underset{t \rightarrow 0}{\rm lim\ }t^{\frac{\alpha}{2}}
\Big\|u(t)\Big\|_{\dot{H}^s_{\mathcal{L}^{1,\infty}}} 
= \underset{t \rightarrow 0}{\rm lim\ }t^{\frac{\alpha}{2}}
\Big\|v(t)\Big\|_{\dot{H}^s_{\mathcal{L}^{1,\infty}}} = 0.
$$
The estimate  \eqref{4.5} is deduced from the inequality \eqref{4.6}. \qed 
\begin{Lem}\label{lem.4.6}
Let $d-1 <s <d$. Then the bilinear operator $B(u, v)(t)$ is continuous from 
$\mathcal{K}_{s,\infty,T} \times \mathcal{K}_{s,\infty,T}$ 
into $\mathcal{K}_{d-1,1,T}$ and we have the inequality
\begin{equation}\label{4.7}
\big\|B(u, v)\big\|_{\mathcal{K}_{d-1,1,T}} \leq
 C\Big\|u\Big\|_{\mathcal{K}_{s,\infty,T}}
\Big\|v\Big\|_{\mathcal{K}_{s,\infty,T}},
\end{equation}
where C is a positive constant and independent of T.
\end{Lem}
\textbf{Proof}.  First, arguing as in Lemma \ref{lem.4.3}, 
we have the following estimates
\begin{gather*}
\big||\xi|^{d-1}\mathcal{F}\big(B(u, v)(t)\big)(\xi)\big| \notag \\
\lesssim \int_0^t|\xi|^{d-1}e^{-(t - \tau)|\xi|^2}|\xi|
\big(|\widehat{u(\tau)}|*|\widehat{v(\tau)}|\big)(\xi)\dif\tau\\
\lesssim \int_0^t|\xi|^{2d-2s}e^{-(t - \tau)|\xi|^2}
\big\|u(\tau)\big\|_{\dot{H}^s_{\mathcal{L}^{1,\infty}}}
\big\|v(\tau)\big\|_{\dot{H}^s_{\mathcal{L}^{1,\infty}}}\dif\tau,
\end{gather*}
this gives the desired result
\begin{gather}
\Big\||\xi|^{d-1}\mathcal{F}\big(B(u, v)(t)\big)(\xi)\Big\|_{L^{\infty,1}_{\xi}}\notag \\
 \lesssim \int_0^t\Big\||\xi|^{2d-2s}e^{-(t - \tau)|\xi|^2}\Big\|_{L^{\infty,1}_{\xi}}
\big\|u(\tau)\big\|_{\dot{H}^s_{\mathcal{L}^{1,\infty}}}
\big\|v(\tau)\big\|_{\dot{H}^s_{\mathcal{L}^{1,\infty}}}\dif\tau  \notag \\
= \int_0^t(t-s)^{s-d}\Big\||\xi|^{2d-2s}e^{-|\xi|^2}\Big\|_{L^{\infty,1}}
\big\|u(\tau)\big\|_{\dot{H}^s_{\mathcal{L}^{1,\infty}}}
\big\|v(\tau)\big\|_{\dot{H}^s_{\mathcal{L}^{1,\infty}}}\dif\tau  \notag \\
 \lesssim \int_0^t (t - \tau)^{\alpha - 1}\tau^{-\alpha}\underset{0 < \eta < t}
{\rm sup}\eta^{\frac{\alpha}{2}}\Big\|u(\eta)
\Big\|_{\dot{H}^s_{\mathcal{L}^{1,\infty}}}\underset{0 < \eta < t}{\rm sup}
\eta^{\frac{\alpha}{2}}\Big\|v(\eta)\Big\|_{\dot{H}^s_{\mathcal{L}^{1,\infty}}}\dif\tau \notag \\
= \underset{0 < \eta < t}{\rm sup}\eta^{\frac{\alpha}{2}}\Big\|u(\eta)
\Big\|_{\dot{H}^s_{\mathcal{L}^{1,\infty}}}\underset{0 < \eta < t}{\rm sup}
\eta^{\frac{\alpha}{2}}\Big\|v(\eta)\Big\|_{\dot{H}^s_{\mathcal{L}^{1,\infty}}}
\int_0^t (t - \tau)^{\alpha- 1}\tau^{-\alpha}\dif\tau \notag \\
\simeq \underset{0 < \eta < t}{\rm sup}\eta^{\frac{\alpha}{2}}
\Big\|u(\eta)\Big\|_{\dot{H}^s_{\mathcal{L}^{1,\infty}}}
\underset{0 < \eta < t}{\rm sup}\eta^{\frac{\alpha}{2}}
\Big\|v(\eta)\Big\|_{\dot{H}^s_{\mathcal{L}^{1,\infty}}}. \label{4.8}
\end{gather}
From \eqref{4.8} it follows \eqref{4.2} since  
$$
\underset{t \rightarrow 0}{\rm lim\ }
\Big\|B(u, v)(t)\Big\|_{\dot{H}^{d-1}_{\mathcal{L}^{1,1}}} 
= \underset{t \rightarrow 0}{\rm lim\ }\Big\||\xi|^{d-1}\mathcal{F}
\big(B(u, v)(t)\big)(\xi)\Big\|_{L^{\infty,1}_{\xi}}  = 0,
$$
whenever
$$
\underset{t \rightarrow 0}{\rm lim\ }t^{\frac{\alpha}{2}}
\Big\|u(t)\Big\|_{\dot{H}^s_{\mathcal{L}^{1,\infty}}} 
= \underset{t \rightarrow 0}{\rm lim\ }t^{\frac{\alpha}{2}}
\Big\|v(t)\Big\|_{\dot{H}^s_{\mathcal{L}^{1,\infty}}} = 0.
$$
The estimate  \eqref{4.7} can be deduced from the inequality \eqref{4.8}. \qed 
\begin{Th}\label{th.4.1} Let $d-1<s<d$ and $1 \leq r <\infty$. 
Then there exists a positive constant $\delta_{s,d}$ such that for 
all $T > 0$ and for all $u_0 \in \dot{H}^{d-1}_{\mathcal{L}^{1,r}}(\mathbb{R}^d),
\ with\ {\rm div}(u_0) = 0$ satisfying
\begin{equation}\label{4.9}
\underset{0 < t < T}{\rm sup}t^{\frac{1}{2}(s+1-d)}
\big\|e^{t\Delta}u_0\big\|_{\dot{H}^s_{\mathcal{L}^1}} \leq \delta_{s,d},
\end{equation}
NSE has a unique mild solution $u \in \mathcal{K}_{s,r,T} \cap L^\infty([0, T]; \dot{H}^{d-1}_{\mathcal{L}^{1,r}})$.\\
In particular, the inequality \eqref{4.9} holds for arbitrary $u_0 \in \dot{H}^{d-1}_{\mathcal{L}^{1,r}}(\mathbb{R}^d)$ 
when $T(u_0)$ is small enough, and there exists a positive 
constant $\sigma_{s,d}$ such that we can take  $T = \infty$ whenever $\Big\|u_0\Big\|_{\dot{H}^{d-1}_{\mathcal{L}^1}} \leq \sigma_{s,d}$.
\end{Th}
\textbf{Proof}. The proof of Theorem \ref{th.4.1} is similar to 
that of Theorem  \ref{th.3.1}. Applying Lemma \ref{lem.4.3} 
and Theorem \ref{th.2.4}, we deduce that there exists a positive 
constant $\delta_{s,d}$ such that
for any $u_0 \in \dot{H}^{d-1}_{\mathcal{L}^{1,r}}(\mathbb{R}^d)$ 
with ${\rm div}(u_0) = 0$ such that
$$
\underset{0 < t < T}{\rm sup}t^{\frac{1}{2}(s+1-d)}
\big\|e^{t\Delta}u_0\big\|_{\dot{H}^s_{\mathcal{L}^{1,\infty}}} 
= \underset{0 < t < T}{\rm sup}t^{\frac{1}{2}(s+1-d)}
\big\|e^{t\Delta}u_0\big\|_{\dot{H}^s_{\mathcal{L}^1}} \leq \delta_{s,d},
$$
the Navier-Stokes equations has a solution $u \in \mathcal{K}_{s,\infty,T}$. 
Applying Lemmas \ref{lem.2.5} and \ref{lem.4.6} we deduce that $u \in L^\infty([0, T]; \dot{H}^{d-1}_{\mathcal{L}^{1,r}})$.  Applying Lemma \ref{lem.4.2}, we get 
$e^{t\Delta}u_0 \in \mathcal{K}_{s,r,T}$. 
From the definition of $\mathcal{K}_{s,r,T}$, we deduce that the left-hand 
side of the inequality \eqref{4.9} converges to $0$ when $T$ 
tends to $0$. Therefore the inequality \eqref{4.9} holds for arbitrary 
$u_0 \in \dot{H}^{d-1}_{\mathcal{L}^{1,r}}(\mathbb{R}^d)$ 
when $T(u_0)$ is small enough. \\
Next, from the inequality \eqref{4.10} with $r=\infty$, we deduce that  
$$
\underset{0 < t < \infty}{\rm sup}t^{\frac{1}{2}(s+1-d)}
\big\|e^{t\Delta}u_0\big\|_{\dot{H}^s_{\mathcal{L}^1}} 
\lesssim  \Big\|u_0\Big\|_{\dot{H}^{d-1}_{\mathcal{L}^1}},
$$
then there exists a positive constant 
$\sigma_{s,d}$ such that  $T = \infty$ and  \eqref{4.9} 
holds whenever $\Big\|u_0\Big\|_{\dot{H}^{d-1}_{\mathcal{L}^1}} 
\leq \sigma_{s,d}$ .\qed 
\begin{Not}
The case $r=\infty$ was studied by Le Jan and Sznitman in \cite{Le Jan  1997}.  
They showed  that NSE are well-posed when  the initial datum belongs to the space  $\dot{H}^{d-1}_{\mathcal{L}^{1,\infty}}$. For $1 \leq r < \infty$ we have 
the following imbedding map  
$$
\dot{H}^{d-1}_{\mathcal{L}^{1,r}}(\mathbb{R}^d)
 \hookrightarrow \dot{H}^{d-1}_{\mathcal{L}^{1,\infty}}(\mathbb{R}^d) 
= \dot{H}^{d-1}_{\mathcal{L}^1}(\mathbb{R}^d).
$$
However, note that for $1 \leq r < \infty$ a function in 
$\dot{H}^{d-1}_{\mathcal{L}^{1,r}}(\mathbb{R}^d)$ 
can be arbitrarily large in the $\dot{H}^{d-1}_{\mathcal{L}^{1,r}}(\mathbb{R}^d)$ 
norm but small in the $\dot{H}^{d-1}_{\mathcal{L}^1}(\mathbb{R}^d)$ norm. 
Theorem \ref{th.4.1} shows the existence of global mild 
solutions in the spaces  
$L^\infty([0, \infty); \dot{H}^{d-1}_{\mathcal{L}^{1,r}}
(\mathbb{R}^d))$  (with $1 \leq r < \infty$)
when the norm of the initial value in the spaces 
$\dot{H}^{d-1}_{\mathcal{L}^1}(\mathbb{R}^d)$ 
is small enough. 
\end{Not}
\vskip 0.5cm 
{\bf Acknowledgments}. This research is funded by 
Vietnam National \linebreak Foundation for Science and Technology 
Development (NAFOSTED) under grant number  101.02-2014.50.

\end{document}